\documentclass[11pt, a4paper]{amsart}

\usepackage{amssymb,amsmath,amsthm}
\usepackage{tikz-cd}
\usetikzlibrary{arrows}
\usepackage{fullpage}
\usepackage{graphicx}
\usepackage{hyperref}
\usepackage{mathtools}
\usepackage[titletoc,title]{appendix}
\usepackage{graphicx}
\usepackage{caption}
\usepackage{subcaption}
\usepackage{blindtext}
\usepackage[titletoc,title]{appendix}
\usepackage[english]{babel}
\usepackage{mathrsfs}
\usepackage{cite}

\newtheorem{theorem}{Theorem}[section]
\newtheorem{lemma}[theorem]{Lemma}
\newtheorem{corollary}[theorem]{Corollary}
\newtheorem{proposition}[theorem]{Proposition}
\theoremstyle{definition}
\newtheorem{definition}[theorem]{Definition}
\newtheorem{example}[theorem]{Example}

\newtheorem{remark}[theorem]{Remark}

\newtheorem*{theorem*}{Theorem}

\newcommand{\isom}{\cong}

\newcommand{\N}{\mathbb{N}}

\DeclareSymbolFontAlphabet{\amsmathbb}{AMSb}

\DeclareMathOperator{\suff}{Suff}
\DeclareMathOperator{\sub}{Sub}
\DeclareMathOperator{\super}{Sup}
\DeclareMathOperator{\pref}{Pref}

\DeclareMathOperator{\ap}{Arith}
\DeclareMathOperator{\ord}{ord}
\DeclareMathOperator{\per}{per}

\DeclareMathOperator{\ideg}{deg^-}
\DeclareMathOperator{\odeg}{deg^+}
\DeclareMathOperator{\bstar}{St^-}
\DeclareMathOperator{\fstar}{St^+}

\DeclareMathOperator{\rk}{rk}

\DeclareMathOperator{\core}{Core}
\DeclareMathOperator{\rcore}{Core^{rel}}

\DeclareMathOperator{\lift}{\mathcal{L}}
\DeclareMathOperator{\olift}{\mathcal{L}_{\text{overlap}}}
\DeclareMathOperator{\slift}{\mathcal{L}_{\text{sub}}}

\DeclareMathOperator{\wmax}{\mathcal{M}}
\DeclareMathOperator{\swmax}{\mathcal{M}_{\text{sub}}}
\DeclareMathOperator{\owmax}{\mathcal{M}_{\text{overlap}}}

\DeclareMathOperator{\subwmax}{\mathcal{SM}}
\DeclareMathOperator{\ssubwmax}{\mathcal{SM}_{\text{sub}}}
\DeclareMathOperator{\osubwmax}{\mathcal{SM}_{\text{overlap}}}

\DeclareMathOperator{\ndir}{\mathsf{ndir}}
\DeclareMathOperator{\odir}{\mathsf{odir}}

\DeclarePairedDelimiter\abs{\lvert}{\rvert}

\makeatletter
\let\oldabs\abs
\def\abs{\@ifstar{\oldabs}{\oldabs*}}

\newcommand{\fakeenv}{} 

\newenvironment{restate}[2]  
{ 
 \renewcommand{\fakeenv}{#2} 
 \theoremstyle{plain} 
 \newtheorem*{\fakeenv}{#1~\ref{#2}} 
 \begin{\fakeenv}
}
{
 \end{\fakeenv}
}

\usepackage{blindtext}
\title{On the intersections of finitely generated subgroups of free groups: reduced rank to full rank}
\author{Marco Linton}

\address{Mathematics institute, Zeeman building, university of Warwick, Coventry, CV4 7AL}

\email{marco.linton@warwick.ac.uk}

\begin{document}
\maketitle

\begin{abstract}
We show that the number of conjugacy classes of intersections $A\cap B^g$, for fixed finitely generated subgroups $A, B<F$ of a free group, is bounded above in terms of the ranks of $A$ and $B$; this confirms an intuition of Walter Neumann. This result was previously known only in the case where $A$ or $B$ is cyclic by the $w$-cycles theorem of Helfer and Wise and, independently, Louder and Wilton. In order to prove our main theorem, we introduce new results regarding the structure of fibre products of finite graphs. We work with a generalised definition of graph immersions so that our results apply to the theory of regular languages. For example, we give a new algorithm to decide non-emptiness of the intersection of two regular languages.
\end{abstract}

\section{Introduction}

Given two combinatorial graph maps $\gamma:\Gamma\to \Delta$ and $\lambda:\Lambda\to \Delta$, we here study their fibre product $\Gamma\times_{\Delta}\Lambda$. This explicitly constructed graph is the pullback in the category of combinatorial graphs. There are two main contexts in which this graph arises: free groups \cite{sta_83} and regular languages \cite{rabin_59}. In the context of free groups, the topological structure of $\Gamma\times_{\Delta}\Lambda$ was intensely investigated in the lead up to the resolution of the strengthened Hanna Neuman conjecture (see in particular \cite{sta_83,neumann_90,tardos_92,dicks_01,khan_02}). In the context of regular languages, the topological and combinatorial structure of $\Gamma\times_{\Delta}\Lambda$ has been studied with applications to free monoids (see \cite{singh_11,giambruno_08}) and to the state complexity of minimal automata (see \cite{birget_92,han_07}). Here we give new results on the structure of the fibre product graph and provide applications in both contexts.

Our first result deals with the simpler case when $\Lambda\isom \sqcup S^1$. A cycle $\lambda:S^1\to \Delta$ is \emph{primitive} if it does not factor properly through any other cycle $S^1\to \Delta$. The \emph{degree} $\deg(\lambda)$ will be the largest degree of a covering map $S^1\to S^1$ that $\lambda$ factors through. Recall that if $X$ is a topological space, then $\beta_i(X)$ is the \emph{$i^{\text{th}}$ Betti number} of $X$. In \cite{helfer_16} and \cite{louder_17} it was shown that if $\gamma:\Gamma\to \Lambda$ is an immersion and $\lambda:S^1\to \Delta$ is primitive, then $\core(\Gamma\times_{\Delta}S^1)$ consists of at most $\beta_1(\Gamma)$ components. This improved upon other bounds given in \cite{wise_05} and resolved the $w$-cycles conjecture, first stated in \cite{wise_03}. If we allow $\lambda$ to be imprimitive, then we would have at most $\beta_1(\Gamma)$ distinct free homotopy classes of cycles instead. Allowing $\lambda$ to also consist of arbitrarily many cycles, we prove the following result.

\begin{restate}{Theorem}{long_cycles_graph_equiv}
Let $\gamma:\Gamma \to \Delta$ and $\lambda:\mathbb{S}\to \Delta$ be graph immersions with $\mathbb{S} = \sqcup S^1$. If $\lambda|_{S^1}$ is primitive for each $S^1\subset \mathbb{S}$, then there are at most $\beta_1(\Gamma)$ many components $S^1\subset \core(\Gamma\times_{\Delta}\mathbb{S})$ such that:
\[
\deg(p_{\mathbb{S}}|_{S^1})\geq 12\cdot \beta_1(\Gamma),
\]
where $p_{\mathbb{S}}:\Gamma\times_{\Delta}\mathbb{S}\to \mathbb{S}$ is the natural projection map.
\end{restate}

From this theorem it follows that there is a unique map $s_{\Gamma}:\mathbb{S}_{\Gamma}\to \Gamma$ where $\mathbb{S}_{\Gamma}$ consists of a disjoint union of at most $\beta_1(\Gamma)$ circles, such that, if $f:S^1\to \Gamma$ is a primitive cycle satisfying $\deg(\gamma\circ f)\geq 12\cdot \beta_1(\Gamma)$, then $f$ factors through $s_{\Gamma}$. In Theorem \ref{long_cycles_algorithm} we give an efficient algorithm to compute the map $s_{\Gamma}$. Theorems \ref{long_cycles_graph_equiv} and \ref{long_cycles_algorithm} are the culmination of several results proved in Section \ref{w_graphs_section} on the full fibre product $\Gamma\times_{\Delta}S^1$. For example, we prove bounds on the number and the sizes of components of $\Gamma\times_{\Delta}S^1$, even when $\gamma$ is not necessarily an immersion.

In Section \ref{main_structure_section} we allow $\Lambda$ to be an arbitrary finite graph. As pointed out in \cite{domaratzki_07} there is a precise multiplicative relation between the sizes of the edge sets of $\Gamma$ and $\Lambda$ and that of their fibre product. If $\Gamma$ and $\Lambda$ represent regular languages or subgroups of a free group, by passing to the core of their fibre product, relative to a collection of \emph{marked vertices}, the subgroup or regular language in question is not changed. If $\Gamma, \Lambda$ and $\Delta$ are all isomorphic to $S^1$, then there is still a multiplicative relation, even after passing to the cores. We show that if we fix the topology of our graphs, this is essentially the only situation in which this relation can occur. The precise statement is as follows.

\begin{restate}{Theorem}{main_structure}
There exists a universal constant $C>0$ such that the following holds. Let $\gamma:\Gamma\to \Delta$ and $\lambda:\Lambda\to \Delta$ be immersions of graphs and let $s:\mathbb{S} = \core(\mathbb{S}_{\Gamma}\times_{\Delta}\mathbb{S}_{\Lambda})\to \core(\Gamma\times_{\Delta}\Lambda)$ be the natural map. Then
\[
\abs{E(\core(\Gamma\times_{\Delta}\Lambda) - s(\mathbb{S}))} \leq C\cdot(\beta_1(\Gamma)\cdot \beta_1(\Lambda))^2\cdot(\abs{E(\Gamma)} + \abs{E(\Lambda)}).
\]
\end{restate}

Finally, we also prove similar bounds to those provided by the $w$-cycles theorem in \cite{helfer_16} and \cite{louder_17} for arbitrary pairs of forwards immersion.

\begin{restate}{Theorem}{full_homotopy_components}
There exists a universal constant $C>0$ such that the following holds. Let $\gamma:\Gamma\to \Delta$ and $\lambda:\Lambda\to \Delta$ be immersions of graphs. The natural map $\theta:\core(\Gamma\times_{\Delta}\Lambda)\to \Delta$ consists of a disjoint union of at most:
\[
C\cdot (\beta_1(\Gamma)\cdot \beta_1(\Lambda))^2
\]
many free homotopy classes of connected graph maps.
\end{restate}

We prove slightly more general versions of these theorems by also taking into account a set of marked vertices and relaxing the condition on the maps. The full statements of Theorems \ref{main_structure} and \ref{full_homotopy_components} also contain explicit upper bounds for their universal constants.

Now that we have stated our main results for the fibre product graph, we begin with some applications.

\subsection{Free groups}

Since Stallings' seminal paper \cite{sta_83}, using finite graphs and graph maps to represent free groups and their finitely generated subgroups has been very fruitful. Therefore, results for graphs and graph maps usually come with complementary results in the theory of free groups. Many examples of such results can be found in \cite{kapovich_02}. We now discuss the consequences of our theorems in the realm of free groups, starting with the following natural application of Theorem \ref{long_cycles_graph_equiv}.

\begin{restate}{Theorem}{finite_index}
Let $F$ be a free group and $A<F$ a finitely generated subgroup. Then there are at most $\rk(A)$ many conjugacy classes of maximal cyclic subgroups $M< F$ such that 
\[
12\cdot\rk(A)\leq[M:A\cap M]<\infty.
\]
\end{restate}

This has an immediate corollary. In \cite{delgado_21}, Delgado, Ventura and Zhakarov introduce the notions of relative order and spectrum. Given a subset $S$ of a group $G$, the \emph{relative order} $\mathcal{O}_S(g)$ of an element $g\in G$, with respect to $S$, is the smallest integer $n\geq 1$ such that $g^n\in S$. The \emph{$S$-spectrum} of $G$ is the set $\mathcal{O}_S(G) = \cup_{g\in G}\mathcal{O}_S(g)$. It was shown in \cite{kapovich_02} that the relative spectrum of a finitely generated subgroup of a free group is bounded. In \cite{delgado_21}, this result is generalised  and the relative spectrum for other classes of groups is also studied. Using Theorem \ref{finite_index}, we may recover the result for free groups and also improve the known bound to one that, in some sense, only depends on the rank of the subgroup. Moreover, by Theorem \ref{long_cycles_algorithm}, we may efficiently compute an upper bound of the $A$-spectrum for finitely generated subgroups $A<F$ of a free group.

\begin{restate}{Corollary}{relative_spectrum}
Let $F$ be a free group and $A<F$ a finitely generated subgroup. Then there exists a set $N\subset \N$ of cardinality bounded above by $13\cdot \rk(A)$, such that:
\[
\mathcal{O}_{A}(F) = \bigcup_{n\in N}\bigcup_{d\mid n} d.
\]
\end{restate}

Our last result of this section is an application of Theorem \ref{full_homotopy_components}. It has been known since the 1950s by work of Howson \cite{howson_54} that the intersection of two finitely generated subgroups of a free group is again finitely generated. If $A$ is a free group, recall that $\bar{\rk}(A) = \min{\{0, \rk(A) - 1\}}$ is the \emph{reduced rank} of $A$. It was conjectured by Hanna Neumann that for any two finitely generated subgroups $A, B<F$ of a free group, the inequality $\bar{\rk}(A\cap B)\leq \bar{\rk}(A)\cdot \bar{\rk}(B)$ held \cite{neumann_57}. A clean graph argument (see \cite{neumann_90}) shows that in fact $\sum_{[A\cap B^g]}\bar{\rk}(A\cap B^g) \leq 2\cdot \bar{\rk}(A)\cdot \bar{\rk}(B)$, where $[A\cap B^g]$ denotes the conjugacy class of the subgroup $A\cap B^g$. In \cite{neumann_90}, Walter Neumann formulated the strengthened Hanna Neumann conjecture, asserting that $\sum_{[A\cap B^g]}\bar{\rk}(A\cap B^g) \leq \bar{\rk}(A)\cdot \bar{\rk}(B)$. In 2011, this long standing conjecture was settled independently by Mineyev \cite{mineyev_12} and Freidmann \cite{friedman_14}.

In the same article in which the strengthened Hanna Neumann conjecture was proposed \cite{neumann_90}, it was also shown that the larger sum $\sum_{[A\cap B^g]}\rk(A\cap B^g)$ is finite. The bound that one can derive from the argument depends on the size of generating sets for $A$ and $B$. This, in turn, depends on a chosen embedding. However, Walter Neumann goes on to suggest that it may be possible to bound this larger sum solely in terms of the ranks of $A$ and $B$.

Since then, most of the literature has focused on the strengthened Hanna Neumann conjecture. However, one improvement to the naive full rank bound in the case where $B = \langle w\rangle$ is cyclic is offered by the $w$-cycles theorem (see \cite{helfer_16,louder_17}). Indeed, it is an immediate consequence that the inequality $\sum_{[A\cap \langle w\rangle^g]}\rk(A\cap \langle w\rangle^g)\leq \rk(A)$ holds. This result is even referred to as a rank 1 version of the Hanna Neumann theorem \cite{helfer_16}. As a consequence of Theorem \ref{full_homotopy_components}, we prove a more general result for arbitrary finitely generated subgroups of a free group, confirming Walter Neumann's intuition.

\begin{restate}{Theorem}{intersection_rank}
There exists a universal constant $C>0$, with the following property. For any free group $F$ and any finitely generated subgroups $A, B<F$, we have:
\[
\sum_{[A\cap B^g]} \rk(A\cap B^g) \leq C\cdot (\rk(A)\cdot \rk(B))^2.
\]
\end{restate}

As demonstrated by Example \ref{lower_bound}, a lower bound for a sharp upper bound for this sum is $\rk(A)\cdot \rk(B)$.

\subsection{Regular languages}

Let $\Delta$ be a rose graph with an edge for each element in some finite set $\Sigma$. Then to each regular language $L\subset \Sigma^*$, there is an associated unique minimal deterministic finite state automaton $\gamma:\Gamma\to \Delta$, such that the language $L(\Gamma)$ is precisely $L$. It is a classical fact that $L(\Gamma\times_{\Delta}\Lambda) = L(\Gamma)\cap L(\Lambda)$ \cite{rabin_59}, thus providing an explicit deterministic quadratic time algorithm to compute the intersection of two regular languages. By considering finite sheeted covers of graphs we see that there are pairs of regular languages that require a quadratic number of states to express their intersection. Hence, this automaton cannot be improved. On the other hand, it is an open question as to whether the non-emptiness of $L(\Gamma)\cap L(\Lambda)$ can be decided in subquadratic time \cite{kara_00}. This decision problem, extending to an input of $k$ DFAs, is known as the $k$-DFA non-emptiness intersection problem ($k$-DFA-NEI). Restricting to unary alphabets, $2$-DFA-NEI is known to be solvable in deterministic linear time  and $3$-DFA-NEI in subquadratic time \cite{oliveira_18}.

If the true complexity of DFA-NEI was strictly better than that of the Rabin-Scott algorithm, then several strong complexity results would be implied. For example, if $2$-DFA-NEI could be solved in subquadratic time, it would imply that $\text{NSPACE}[n + o(n)]\subset \text{DTIME}[2^{(1-\epsilon)\cdot n}]$ for some $\epsilon>0$, SAT would be solvable in subexponential time and the strong exponential time hypothesis would be contradicted \cite{oliveira_20}. If $k$-DFA-NEI could be solved in $O(n^{o(k)})$ time, we would also have $\text{NSPACE}[n]\subset \text{DTIME}[2^{o(n)}]$ \cite{oliveira_20} and $\text{NL}\neq \text{P}$ \cite{wehar_14}.

We present a different approach to the non-emptiness intersection problem based on Theorem \ref{main_structure}. If we fix the topology of the automata, we obtain complexity bounds that depend linearly on the input size. Define $\partial\Gamma\subset V(\Gamma)$ to be the vertices with indegree or outdegree equal to $0$.

\begin{restate}{Theorem}{intersection_emptiness}
\label{intersection_emptiness_main}
For some $q>0$, there is an algorithm that, given as input two DFAs $\gamma:\Gamma\to \Delta$ and $\lambda:\Gamma\to \Delta$ such that $\beta_1(\Gamma), \beta_1(\Lambda), \abs{\partial\Gamma}, \abs{\partial\Lambda}\leq m$ and $\abs{E(\Gamma)}, \abs{E(\Lambda)}\leq n$, decides if $L(\Gamma)\cap L(\Lambda)$ is non-empty in
\[
O(m^q\cdot n)
\]
time.
\end{restate}

Note that if $\Delta$ contains a single edge, then for every forwards graph immersion $\gamma:\Gamma\to \Delta$ we have $\beta_1(\Gamma) + \abs{\partial\Gamma}\leq 2$. Thus, this can be seen as a direct generalisation of the unary case.

\subsection{Organisation of the paper}

Our main technical tools will come from new results on the combinatorics of words. In Section \ref{words_section}, we recall some well known definitions and results and introduce $w$-factorisations and word overlaps. We then show uniqueness of $w$-factorisations and describe all overlaps for any given pair of words in terms of arithmetic progressions. The main result of this section is Theorem \ref{cyclic_shift}, upon which most of Section \ref{w_graphs_section} relies. In Section \ref{graphs_section} we recall all the necessary definitions for graphs and the fibre product graph. Section \ref{w_graphs_section} introduces $w$-graphs and culminates in the proof of Theorem \ref{long_cycles_graph_equiv}. In Section \ref{main_structure_section}, we decompose the fibre product into two subgraphs. Using this decomposition and combining our word overlap and $w$-graph results, we then obtain Theorems \ref{main_structure} and \ref{full_homotopy_components}. Since the graphs we work with are not Stallings graphs, in Section \ref{free_groups_section}, we explain an alternative construction that allows us to convert our results to statements about free groups. Finally, Section \ref{algorithm_section} is dedicated to our algorithmic results.

\subsection*{Acknowledgments} We would like to thank Saul Schleimer for many helpful discussions and comments on previous versions of this article.

\section{Some combinatorics on words}
\label{words_section}

\subsection{Words and the free monoid}

Let $\Sigma$ be a set. We denote by $\Sigma^*$ the free monoid generated by $\Sigma$ under the concatenation operation and with neutral element the empty string $\epsilon$. We call elements of $\Sigma$ \emph{letters} and elements of $\Sigma^*$ \emph{words}. 

We call an equality of the form $w = w_0\cdot w_1\cdot...\cdot w_n\in \Sigma^*$ a \emph{factorisation} of $w$ and each word $w_i$ a \emph{subword} of $w$. Now let $w\in \Sigma^*$ and 
\[
w = w_0\cdot w_1\cdot...\cdot w_n
\]
be the unique factorisation of $w$ with $w_i\in \Sigma$. Then we denote by $\abs{w} = n + 1$ the \emph{length} of $w$ and by $w[i] = w_i$ the $i+1^{\text{th}}$ letter of $w$. For any $0\leq i< j\leq n+1$, we denote by $w[i:j] = w_i\cdot w_{i+1}\cdot...\cdot w_{j-1}$ the subword of length $j - i$ starting at $w_i$. We denote by $w[:i] = w[0:i]$ the \emph{prefix} of length $i$ and by $w[j:] = w[j:n+1]$ the \emph{suffix} of length $n - j + 1$. We extend this all to negative indices by defining $w[-j:-i] = w[|w| - j:|w| - i]$. A \emph{cyclic conjugate} of $w$ is a word of the form $w[-i:]\cdot w[:i]$ for some $i\leq \abs{w}$. We say that two words $v, w$ are \emph{cyclically conjugate} if $v$ is a cyclic conjugate of $w$. A \emph{conjugating word} of two cyclically conjugate words $v, w$ is a word $u$ such that $u\cdot v = w\cdot u$.

A word $w$ is \emph{primitive} if there is no subword $u$ such that $w = u^n$ for some $n>1$, \emph{imprimitive} otherwise. A \emph{period} $p$ of a word $w$ is an integer $p<\abs{w}$ such that $w[i] = w[i+p]$ for all $i\leq |w|-p$. Denote by $\per(w)$ the smallest period of $w$. Let $w[:\per(w)] = u$, then $w = u^n\cdot u'$ with $u'$ a proper prefix of $u$. We will call $n = \ord(w)$ the \emph{order} of $w$. Denote by $\ap(p, q, r)$ the \emph{arithmetic progression} $\{i\cdot q + r | 0\leq i\leq p\}$.

\subsection{Periodicity and $w$-factorisations}

Before introducing $w$-factorisations, we will state some results on periodic words, starting with the well known Fine and Wilf periodicity theorem \cite{fine_65}.

\begin{theorem}
\label{periodicity}
If $w\in \Sigma^*$ has periods $p$ and $q$ and $|w|\geq p + q - \gcd(p, q)$, then $w$ has period $\gcd(p, q)$.
\end{theorem}

\begin{corollary}
\label{square_period}
Let $w\in \Sigma^*$. If $w^2$ has period $p<\abs{w}$, then $w$ is imprimitive.
\end{corollary}

The following are Lemmas 1 and 3 of \cite{lyndon_62} respectively.

\begin{lemma}
\label{key_power}
Let $u, v, w\in \Sigma^*$ such that $u\cdot w = w\cdot v$. Then either $w$ is a proper prefix of $u$, or $w$ has prefix $u$ and period $|u|$.
\end{lemma}

\begin{lemma}
\label{commuting}
Let $v, w\in \Sigma^*$ such that $v\cdot w = w\cdot v$. Then $v$ and $w$ are powers of a common word.
\end{lemma}

A key tool in this paper will be the notion of $w$-factorisations:

\begin{definition}
Suppose $w\in \Sigma^*$ is a primitive word. A \emph{left $w$-factorisation} of a word $x\in \Sigma^*$ is a factorisation of the form $x = v\cdot y$ where $v$ is a subword of $w^n$ for some $n\geq 1$. A left $w$-factorisation $x = v\cdot y$ is \emph{maximal} if there is no other left $w$-factorisation $x = v'\cdot y'$ with $\abs{v'}>\abs{v}$. A \emph{(maximal) right $w$-factorisation} is defined analogously.
\end{definition}

\begin{lemma}
\label{unique_factorisation}
Suppose $w\in \Sigma^*$ is primitive. Then for any word $x\in \Sigma^*$, we have:
\begin{enumerate}
\item left and right $w$-factorisations of $x$ of length $k\geq |w|$ are unique,
\item if $x$ has a maximal left $w$-factorisation of length $k_1$ and a maximal right $w$-factorisation of length $k_2$, then either $k_1 = k_2 = \abs{w}$, or $k_1 + k_2<|x| + |w|$.
\end{enumerate}
\end{lemma}

\begin{proof}
If the first statement does not hold, then up to replacing $w$ by a cyclic conjugate of itself, we get $x = w\cdot y = v\cdot u \cdot y$ where $v$ is a proper suffix of $w$ and $u$ is a proper prefix of $w$. But then $w = u\cdot v = v\cdot u$ so by Lemma \ref{commuting} $w$ is not primitive. The second statement follows from the same argument.
\end{proof}

\begin{remark}
Lemma \ref{unique_factorisation} cannot be extended for any left or right factorisations of length less than $w$ for the simple fact that $w$ can have many repeated subwords. For example, if $w = u^n\cdot u'$ for some $n\geq 1$ and $u'$ a proper non-trivial prefix of $u$, then $x = u^{n-k}\cdot u'$ has at least $k+1$ different $w$-factorisations.
\end{remark}

We conclude this section with a couple of algorithmic results. The following theorem is the main result in \cite{knuth_74}.

\begin{theorem}
\label{linear_pattern}
The set of periods $p$ of a word $w\in \Sigma^*$ can be computed in linear time.
\end{theorem}

\begin{proposition}
\label{cyclic_conjugation_algorithm}
Let $u, v\in \Sigma^*$. We can check if $u$ and $v$ are cyclically conjugate in $O(\abs{u} + \abs{v})$ time. Moreover, the algorithm returns a conjugating word $w$.
\end{proposition}

\begin{proof}
The words $u$ and $v$ are cyclically conjugate if and only if $\abs{u} = \abs{v}$ and $v$ is a subword of $u^2$. Since checking this can be done in linear time by the Knuth-Morris-Pratt algorithm \cite{knuth_74}, checking if $u$ and $v$ are cyclically conjugate can be done in $O(\abs{u} + \abs{v})$ time. Since the algorithm also outputs an index $i$ such that $u^2[i:i+\abs{v}] = v$, we also get the conjugating word $w = u[:i]$.
\end{proof}

\subsection{Overlapping words}
\label{word_overlaps}

In this section we will review the idea of overlaps of two given words $v, w\in \Sigma^*$. 

\begin{definition}
Let $v, w\in \Sigma^*$ be two words. There are three distinct types of overlaps: \emph{subword}, \emph{prefix} and \emph{suffix} overlap. These are described in the table below.
\begin{align*}
v[i:j] &= w && \text{Subword overlap}\\
v[:i] &= w[-i:] && \text{Prefix overlap}\\
v[-i:] &= w[:i] && \text{Suffix overlap}
\end{align*}
We will use the notation $\sub(v, w)$ to denote the set of \emph{subword indices of $w$ in $v$}: that is, the set of integers $0\leq i\leq \abs{v} - \abs{w}$ such that $v[i:i+|w|] = w$. We will write $\pref(v, w)$ to denote the set of \emph{prefix indices of $w$ in $v$}: that is, the set of integers $0 < i\leq \min{\{\abs{v}, \abs{w}\}}$ such that $v[:i] = w[-i:]$. Finally, we write $\suff(v, w)$ to denote the set of \emph{suffix indices of $w$ in $v$}: that is, the set of integers $\max{\{0, \abs{v} - \abs{w}\}}\leq i <\abs{v}$ such that $v[i:] = w[:\abs{v} - i]$.
\end{definition}

\begin{definition}
We say $w$ \emph{crosses} $v$ at $j$ if there is an integer $i$ with $i\leq  j< i + |w|$ such that
\[
v[i : i + |w|] = w.
\]
We will call the index $i$ a \emph{crossing index of $w$ in $v$ at $j$}. Denote by $\sub(v, w, j)$ the crossing indices of $w$ in $v$ at $j$.
\end{definition}

The word $w$ crosses $v$ at $j$ if $w$ appears as a subword of $v$ containing the letter $v[j]$. Or in other words, $w$ crosses $v$ at $j$ if $v = y\cdot w\cdot z$ with $i = |y| + 1$ a crossing index of $w$ in $v$ at $j$. This definition and the following lemma have been used extensively in the pattern matching literature; see \cite{lifshits_07} and the references therein.

\begin{lemma}
\label{crossing}
Let $w, v\in \Sigma^*$ and let $j\leq |v|$ be an integer. There exists some triple $p, q, r\leq |v|$ such that:
\[
\sub(v, w, j) = \ap(p, q, r)
\]
with $q\geq \per(w)$ and $p\leq \ord(w)$. Furthermore, if $p\geq 2$, then $q = \per(w)$.
\end{lemma}

\begin{proof}
The first part of the statement is Lemma 1 in \cite{lifshits_07}. If $p\geq 2$, then we must have $\per(w) | q$ by Theorem \ref{periodicity}. Set $d = \frac{q}{\per(w)}$. If we have
\[
v = y\cdot w[:q]^{p}\cdot w[q:]\cdot z
\]
then we also have
\[
v = y\cdot w[:\per(w)]^{p\cdot d}\cdot w[q:]\cdot z.
\]
Thus $d = 1$ as otherwise $\ap(p, q, r)\subsetneq \sub(v, w, j)$.
\end{proof}

The following results are direct applications of Lemma \ref{crossing}.

\begin{lemma}
\label{subwords}
Let $v, w\in \Sigma^*$. There is a sequence of integers $j_1, ..., j_k$ such that $j_{i+1} - j_i\geq\abs{w}$ for all $i<k$ such that:
\[
\sub(v, w) = \bigcup_{i=1}^k \sub(v, w, j_i).
\]

\end{lemma}

\begin{corollary}
Let $v, w\in \Sigma^*$. The set $\sub(v, w)$ can be decomposed:
\[
\sub(v, w) = \bigcup_{i=1}^k \ap(p_i, q_i, r_i)
\]
where $k\leq \left\lfloor\frac{|v|}{|w|}\right\rfloor$ and $p_i\cdot q_i\leq \abs{w}$ for all $i\leq k$.

\end{corollary}

We now move onto classifying prefix and suffix overlaps in a similar way. Since $\pref(w, v) = \{\abs{v} - i \mid i\in \suff(v, w)\}$, it is only necessary to treat the $\suff(v, w)$ case.

\begin{lemma}
\label{suffix_ap}
Let $v, w\in \Sigma^*$. If $i < j\in \suff(v, w)$, then 
\[
\ap\left(\left\lfloor\frac{\abs{v} - i}{j - i}\right\rfloor, j - i, i\right)\subset \suff(v, w).
\]
\end{lemma}

\begin{proof}
If $i < j\in \suff(v, w)$, then $v[i:] = w[: \abs{v} - i]$ and $v[j:] = w[:\abs{v} - j]$. In particular:
\[
v[i:j]\cdot w[: \abs{v} - i] = w[: \abs{v} - i]\cdot w[\abs{v} - j:\abs{v} - i].
\]
By Lemma \ref{key_power}, $v[i:]$ and $w[:\abs{v} - i]$ have period $j - i$ and so $\ap\left(\left\lfloor\frac{\abs{v} - i}{i - j}\right\rfloor, j - i, i\right)\subset \suff(v, w)$ as claimed.
\end{proof}

\begin{proposition}
\label{suffix_overlap}
Let $v, w\in \Sigma^*$. The set of suffix indices of $w$ in $v$ can be decomposed as a union of arithmetic progressions:
\[
\suff(v, w) = \bigcup_{i=1}^k \ap(p_i, q_i, r_i)
\]
such that, for all $n\leq k$:
\begin{enumerate}
\item\label{itm:cover} if $j\in \suff(v, w)$ and $j\leq |v| - q_n$, then $j\in \bigcup_{i=1}^n\ap(p_i, q_i, r_i)$,
\item\label{itm:start} $r_n \geq |v| - q_{n-1}$,
\item\label{itm:end} $p_n\cdot q_n + r_n \geq |v| - q_n$,
\item\label{itm:growth} $q_n\leq \frac{q_{n-1}}{2}$,
\item\label{itm:prim} $v[r_n : r_n + q_n]$ is primitive.
\end{enumerate}
Furthermore, $k\leq \min{\{\left\lfloor \log_2{|v|}\right\rfloor, \left\lfloor \log_2{|w|}\right\rfloor\}}$
\end{proposition}

\begin{proof}
We will inductively define arithmetic progressions $\ap(p_n, q_n, r_n)\subset \suff(v, w)$ satisfying (\ref{itm:cover}) to (\ref{itm:prim}) in the statement. Once we have done this, (\ref{itm:cover}) will imply the decomposition and (\ref{itm:growth}) will imply the bound on $k$.

For the base case we will simply set $r_1 = \min{\{\suff(v, w)\}}$, $p_1 = 0$ and $q_1 = \abs{v} - r_1$. Now suppose $n\geq 2$. Let $j\in \suff(v, w)$ be the smallest index not contained in $\bigcup_{i=1}^{n-1}\ap(p_i, q_i, r_i)$. Let $l\in \suff(v, w)\cup \{\abs{v}\}$ minimise $\abs{l - j}$. Inductively by (\ref{itm:cover}), we have $j\geq \abs{v} - q_{n-1}$ and by (\ref{itm:end}) we have $p_{n-1}\cdot q_{n-1} + r_{n-1}\geq \abs{v} - q_{n-1}$, hence $\abs{l - j}\leq \frac{q_{n-1}}{2}$. So if we set $r_n = j$, $q_n = \abs{l - j}$ and $p_n = \left\lfloor\frac{\abs{v} - r_n}{q_n}\right\rfloor$, then $\ap(p_n, q_n, r_n)\subset \suff(v, w)$ by Lemma \ref{suffix_ap}. This verifies the last four properties, so we turn our attention to (\ref{itm:cover}). 

If $l = |v|$, then we are done. Now assume that $l\neq \abs{v}$ and hence $p_{n}\geq 1$. Suppose there is some index $m\in \suff(v, w)$ such that $|v| - q_{n-1}< m \leq |v| - q_{n}$ and $m\notin \cup_{i=1}^n\ap(p_{i}, q_{i}, r_{i})$. By our choice of $j$, we must have $m> j$. Furthermore, $\abs{m - j}\geq q_n$ and $q_n\nmid \abs{m - j}$. So in particular, there is some integer $i< p_n$ such that $0<m - (i\cdot q_n + r_n)< q_n$. We may assume that $m$ was chosen to minimise $\abs{m - i\cdot p_n + r_n}$ over all valid choices of $i$. Then by Lemma \ref{suffix_ap}, $v[i\cdot q_n + r_n:]$ has period $m - (i\cdot q_n + r_n)$. Now by Theorem \ref{periodicity} we arrive at a contradiction to the minimality of $m$. Hence no such $m$ exists and (\ref{itm:cover}) holds.
\end{proof}

We close this section with a result which will play a key role in our analysis of $w$-graphs in Section \ref{w_graphs_section}.

\begin{theorem}
\label{cyclic_shift}
Let $w, x_1, ..., x_n\in \Sigma^*$ be words with $w$ primitive. There is some cyclic shift $w_0$ of $w$ such that
\[
|\sub(w_0^2, x_i, |w_0| - 1)| \leq 2
\]
for all $1\leq i\leq n$.
\end{theorem}

In other words, Theorem \ref{cyclic_shift} tells us that there is some cyclic shift $w_0$ of $w$ such that each $x_i$ crosses the midpoint of $w_0^2$ at most twice.

\begin{proof}[Proof of Theorem \ref{cyclic_shift}]
Suppose there is some $x_i$ such that $\abs{x_i}\geq \abs{w}$. By Lemma \ref{unique_factorisation}, if $\sub(w_0^2, x_i, \abs{w} - 1)\neq \emptyset$ for some cyclic shift $w_0$ of $w$, then $x_i$ has a unique $w$-factorisation. Hence, $\abs{\sub(w_0^2, x_i, \abs{w} - 1)}\leq 1$ for all cyclic shifts $w_0$ of $w$. By Lemma \ref{crossing}, if $\ord(x_i)\leq 1$, then we also have $|\sub(w_0^2, x_i, |w| - 1)|\leq 2$ for all $i$. Hence there is no loss in generality assuming that $\abs{x_i}<\abs{w}$ and $\ord(x_i)\geq 2$ for all $i$.

Denote by $w[[j]] = w[j:]\cdot w[:j]$ where $j$ is taken mod $\abs{w}$. For each $i$, let $y_i, y_i', y_i''\in \Sigma^*$ and $c_i\geq 2$ such that 
\begin{align*}
y_i &= x_i[:\per(x_i)]\\
x_i &= y_i^{c_i}\cdot y_i'\\
y_i &= y_i'\cdot y_i''
\end{align*}

We may assume that the $x_i$ are ordered so that $|y_i|\geq |y_{i+1}|$. To simplify the argument we will consider an extra word $x_0 = w$. We will find by induction a sequence of integers $i_0, i_1, ..., i_n$ such that for all $k\leq n$ the following is satisfied:
\begin{enumerate}
\item\label{itm:bounds} if $k\geq 1$, then $i_{k-1} - 2\cdot \abs{y_{k-1}}\leq i_{k} - 2\cdot \abs{y_{k}}$ and $i_{k-1}\leq i_k$,
\item\label{itm:cyclic_sub} $\abs{\sub(w^2[[i_k - l]], x_k, \abs{w} - 1)}\leq 2$ for all $l\leq 2\cdot \abs{y_k}$,
\item\label{itm:smaller_k} there is some integer $m(k)\leq k$ such that $i_{m(k)} = ... = i_{k}$ and $y_{m(k)}^2$ is a suffix of $w[[i_{m(k)}]]$.
\end{enumerate}
Then we will take $w_0 = w[[i_n]]$ and the theorem will be proved by (\ref{itm:bounds}) and (\ref{itm:cyclic_sub}).

By setting $i_0 = 2\cdot \abs{w}$ and $m(0) = 0$, the base case holds trivially. So now for the inductive step, let $k\geq 1$. If
\[
\abs{\sub(w^2[[i_{m(k - 1)} - l]], x_k, \abs{w})}\leq 2
\]
for all $l\leq 2\cdot \abs{y_{m(k-1)}}$, then set $i_k = i_{m(k)} = i_{m(k-1)}$ and $m(k) = m(k-1)$ and all three conditions are satisfied. So now suppose that:
\[
|\sub(w^2[[i_{m(k-1)} - l]], x_k, |w|)|> 2
\]
for some $l\leq2\cdot |y_{m(k-1)}|$. Then $y_k^2$ is a subword of $w[[i_{m(k-1)}]][ - 2\cdot \abs{y_{m(k-1)}}:]$ by Lemma \ref{crossing} and $y_k\neq y_{m(k-1)}$. So let $i, j$ be integers such that
\[
i_{m(k-1)} - 2\cdot \abs{y_{m(k-1)}} \leq i< j\leq i_{m(k-1)}
\]
and $w[[i_{m(k-1)}]][i:j]$ has the largest $y_k$-factorisation among all such choices of integers. If $i = i_{m(k-1)} - 2\cdot \abs{y_{m(k-1)}}$ and $j = i_{m(k-1)}$, then, either $y_{k} = y_{m(k-1)}$, or $y_{m(k-1)}$ is imprimitive by Corollary \ref{square_period}. But this is a contradiction.

So now suppose first that $i> i_{m(k-1)} - 2\cdot \abs{y_{m(k-1)}}$. Let $i_k>i$ be the smallest integer such that $y_k^2$ is a suffix of $w[[i_k]]$, but not $y_k^3$. This exists and is smaller than $j$ since $j-i\geq 2\cdot \abs{y_k}$. Thus we have
\begin{align*}
i_{m(k-1)} - 2\cdot \abs{y_{m(k-1)}}&< i \leq i_k - 2\cdot \abs{y_k}\\
i_k &\leq j \leq i_{m(k-1)}
\end{align*}
hence (\ref{itm:bounds}) holds. By setting $m(k) = k$, (\ref{itm:smaller_k}) holds. If $c\in \sub(w^2[[i_k]], x_k, \abs{w}-1)$, then $w[[i_k]][c:]$ has a $y_k$-factorisation beginning with $y_k$. But then by Lemma \ref{unique_factorisation} it follows that either $c = \abs{w} - 2\cdot \abs{y_k}$ or $c \geq \abs{w} - \abs{y_k}$. So by Lemma \ref{crossing}, (\ref{itm:cyclic_sub}) holds.

If $i = i_{m(k-1)} - 2\cdot \abs{y_{m(k-1)}}$, then $j< i_{m(k-1)}$. So now we let $i_k\leq j$ be the largest integer such that $y_k^2$ is a suffix of $w[[i_k]]$, but $y_k$ is not a prefix. As before, (\ref{itm:bounds}) is satisfied. Similarly, setting $m(k) = k$, (\ref{itm:smaller_k}) is also satisfied. If $c\in \sub(w^2[[i_k]], x_1, \abs{w}-1)$, then by Lemma \ref{crossing}, $w[[i_k]][c:]$ has a $y_k$-factorisation. By Lemma \ref{unique_factorisation}, we must have $c = \abs{w} - d\cdot \abs{y_k}$ for some $d\geq 1$ or $c>\abs{w} - \abs{y_k}$. Since $w[[i_k]]$ does not have $y_k$ as a prefix, $\sub(w^2[[i_k]], x_k, \abs{w} - 1)$ contains at most one value of the first type. By Lemma \ref{crossing} it contains at most one value of the second type. So (\ref{itm:cyclic_sub}) holds and the theorem is proved.
\end{proof}

\section{Graphs}
\label{graphs_section}

A directed graph is a tuple $\Gamma = (V, E, o, t)$ where $V$ is the set of \emph{vertices}, $E$ is the set of \emph{edges} and $o, t:E\to V$ are the \emph{origin} and \emph{target} maps respectively. We always assume our graphs are directed and so will refer to them simply as graphs. If $S\subset V$ or $S\subset E$, then we denote by $\Gamma - S$ the maximal subgraph of $\Gamma$ that does not contain any vertices or edges in $S$. Define the \emph{forwards star} of a vertex $v$ in a graph $\Gamma$ to be the set $\fstar(v, \Gamma) = \{e\in E(\Gamma) \mid o(e) = v\}$. We also define the set $\bstar(v, \Gamma) = \{e\in E(\Gamma) \mid t(e) = v \}$ which we will call the \emph{backwards star}. We then define the \emph{indegree} of $v$ to be $\ideg(v) = \abs{\bstar(v, \Gamma)}$, the \emph{outdegree} of $v$ to be $\odeg(v) = \abs{\fstar(v, \Gamma)}$ and the \emph{degree} of $v$ to be $\deg(v) = \ideg(v) + \odeg(v)$. 

With applications in mind, it is useful to allow an auxiliary subset $V^m(\Gamma)\subset V(\Gamma)$ of \emph{marked vertices}. We further decompose this set $V^m(\Gamma) = V^i(\Gamma)\cup V^f(\Gamma)$ into \emph{initial} and \emph{final vertices}.

\begin{definition}
Denote by $\partial \Gamma\subset V(\Gamma)$ the vertices $v$ such that $\ideg(v) = 0$ or $\odeg(v) = 0$.
\end{definition}

Note that this is a different definition to the usual boundary of a graph.

\begin{definition}
Denote by $\bar{V}(\Gamma)\subset V(\Gamma)$ the vertices $v$ such that $\ideg(v) \neq 1$ or $\odeg(v)\neq 1$.
\end{definition}

A \emph{graph map} $\gamma:\Gamma\to \Delta$ consists of a pair of maps $\gamma_V:V(\Gamma)\to V(\Delta)$ and $\gamma_E:E(\Gamma)\to E(\Delta)$ such that $\gamma_V\circ o = o\circ \gamma_E$ and $\gamma_V\circ t = t\circ \gamma_E$. If $\gamma:\Gamma\to \Delta$ is a graph map of graphs with marked vertices, we require that $V^i(\Gamma)$ maps to $V^i(\Delta)$ and $V^f(\Gamma)$ maps to $V^f(\Delta)$. We say a graph map $\gamma:\Gamma\to\Delta$ is a \emph{forwards immersion} if the induced map $\gamma_v^+: \fstar(v, \Gamma)\to \fstar(\gamma(v), \Delta)$ is injective for each $v\in V(\Gamma)$. We say $f$ is an \emph{immersion} if $\gamma_v^+$ and $\gamma_v^-$ are injective for all $v\in V(\Gamma)$.

A \emph{segment} is a graph consisting of vertices $v_0, ..., v_{k}$ and edges $e_0, ..., e_{k-1}$ such that $o(e_i) = v_i$, $t(e_i) = v_{i+1}$ for all $i< k$. A \emph{path} is a graph map $f:I\to \Gamma$ where $I$ is a segment. The \emph{length} of a path $f$ will be defined as $|f| = |E(I)| = k$. Denote by 
\begin{align*}
f(i) &= v_i,& f[i] &= e_i,\\
o(f) &= v_0,&  t(f) &= v_k.
\end{align*}
Paths may also be considered as words over $E(\Gamma)$ and thus we will be able to apply the same notation and results from Section \ref{words_section} to paths. If $f$ and $g$ are two paths with $t(f) = o(g)$, then we will write $f*g$ for their concatenation. A path $f$ is \emph{primitive} if there is no other path $w:I\to \Gamma$ such that $f = w^n$ for some $n\geq 2$. A path $f$ is \emph{simple} if $f(i) \neq f(j)$ for all $i, j<\abs{f}$. A path $f:I\to \Gamma$ is a \emph{loop} if $o(f) = t(f)$. We say that a path $f:I\to \Gamma$ is \emph{accepting} if $o(f) \in V^i(\Gamma)$ and $t(f)\in V^f(\Gamma)$. A path $f:I\to \Gamma$ is \emph{admissible} if it is accepting or if it is a loop. 

We say a path $f:I\to \Gamma$ \emph{connects} $o(f)$ with $t(f)\in V(\Gamma)$. Connectedness generates a reflexive and transitive binary relation on $V(\Gamma)$. A graph is \emph{connected} it has a single equivalence class under the symmetric closure of this binary relation. The \emph{connected components} are the maximal connected subgraphs.

A \emph{circle} is a graph consisting of vertices $v_0, ..., v_{k-1}$ and edges $e_0, ..., e_{k-1}$ such that $o(e_i) = v_i$ and $t(e_i) = v_{i+1}$ for all $i< k$ and where the indices are taken mod $k$. A \emph{cycle} is a graph map $f:S^1\to \Gamma$ where $S^1$ is a circle. Recall that $\deg(f)$ is the largest degree of a covering map $S^1\to S^1$ that $f$ factors through. If $f:I\to \Gamma$ is a path, we write $\hat{f}:S^1\to \Gamma$ for the cycle obtained from $f$ by identifying the vertices $\partial I$. The union of images of cycles in $\Gamma$ is \emph{the core} of $\Gamma$; we denote this by $\core(\Gamma)$. The \emph{relative core} $\rcore(\Gamma)$ is the union of the images of all admissible paths.

\begin{definition}
Denote by $\bar{E}(\Gamma)$ the collection of paths $e:I\to \Gamma$ and cycles $f:S^1\to \Gamma$ such that:
\begin{enumerate}
\item $e_E$ and $f_E$ are injective,
\item $e(I)\cap \bar{V}(\Gamma) = \{o(e), t(e)\}$ and $f(S^1)\cap \bar{V}(\Gamma) = \emptyset$.
\end{enumerate}
We denote by $\bar{E}_c(\Gamma)\subset \bar{E}(\Gamma)$ the subset corresponding to all the cycles. 
\end{definition}

\begin{figure}
\centering
\vspace*{8pt}
\includegraphics[scale=0.8]{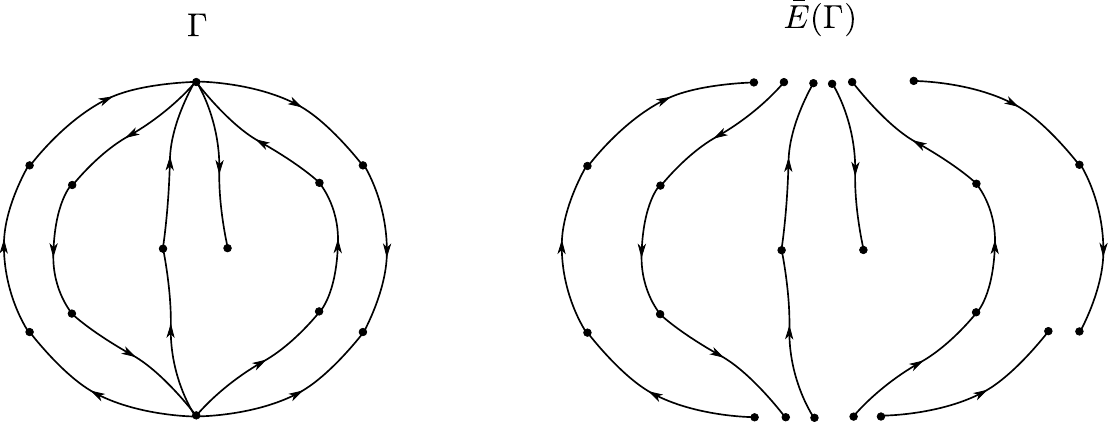}
\vspace{8pt}
\caption{A graph $\Gamma$ and the segments $\bar{E}(\Gamma)$.}
\label{bar_edges}
\end{figure}

See Figure \ref{bar_edges} for an example of the elements $\bar{E}(\Gamma)$. Throughout, we may consider the elements of $\bar{E}(\Gamma)$ as paths or cycles, or as subgraphs of $\Gamma$.

We define an alternative version of our origin and target maps for graphs $\Gamma$ such that $\bar{E}_c(\Gamma) = \emptyset$. If $e\in E(\Gamma)$ and $\bar{e}\in \bar{E}(\Gamma)$ is the unique element such that $e\subset \bar{e}(I)$, then we do this as follows:
\[
\bar{o}:E(\Gamma)\to \bar{E}(\Gamma)\times \N,
\]
\[
\bar{o}(e) = (\bar{e}, i),
\]
where $o(e) = \bar{e}(i)$. Then the map
\[
\bar{t}:E(\Gamma)\to \bar{E}(\Gamma)\times \N
\]
is defined similarly. We also extend these maps to non-trivial paths.

\begin{definition}
Suppose $\gamma:\Gamma\to \Delta$ is a graph map and $w:I\to \Delta$ is a non-trivial path. Then we will denote by $\lift(\gamma, w)$ the lifts of $w$ to $\Gamma$.
\end{definition}

\begin{definition}
Suppose $\gamma:\Gamma\to \Delta$ is a graph map and $w:I\to \Delta$ is a path. Define $\slift(\gamma, w)\subset \lift(\gamma, w)$ to be the lifts that do not contain any vertex $v\in \bar{V}(\Gamma)$ in their image. Define $\olift(\gamma, w)\subset \lift(\gamma, w)$ to be the lifts which do contain some $v\in \bar{V}(\Gamma)$ in their image.
\end{definition}

In Section \ref{counting_lifts}, we characterise these subsets for the fibre product graph.

The following lemma is a rephrasing of a standard fact about graphs using our notation.

\begin{lemma}
\label{betti_bound}
Let $\Gamma$ be a connected graph, then:
\begin{align*}
\abs{\bar{V}(\Gamma)} &\leq \max{\{1, -2\cdot\chi(\Gamma) + 2\cdot\abs{\partial\Gamma}\}},\\
\abs{\bar{E}(\Gamma) - \bar{E}_c(\Gamma)} &\leq \max{\{0, -3\cdot \chi(\Gamma) + 2\cdot \abs{\partial\Gamma}\}}.
\end{align*}

\end{lemma}

\subsection{The fibre product graph}

Let $\gamma:\Gamma\to \Delta$ and $\lambda:\Lambda\to \Delta$ be two graph maps. Define the \emph{fibre product graph} of $\Gamma$ and $\Lambda$ as the graph $\Gamma\times_{\Delta}\Lambda$ with  vertex set:
\[
V(\Gamma\times_{\Delta}\Lambda) = \{(u, v)\in V(\Gamma)\times V(\Lambda) \mid \gamma_V(u) = \lambda_V(v)\}
\]
and edge set:
\[
E(\Gamma\times_{\Delta}\Lambda) = \{(e, f)\in E(\Gamma)\times E(\Lambda) \mid \gamma_E(e) = \lambda_E(f)\}.
\]
The origin and target maps $o, t:E(\Gamma\times_{\Delta}\Lambda)\to V(\Gamma\times_{\Delta}\Lambda)$ will then be given by $o(e, f) = (o(e), o(f))$ and $t(e, f) = (t(e), t(f))$. We will denote by $p_{\Gamma}, p_{\Lambda}:\Gamma\times_{\Delta}\Lambda\to \Gamma, \Lambda$ the natural projection maps. If our graphs have marked vertices, then we will define the marked vertices of $\Gamma\times_{\Delta}\Lambda$ to be 
\begin{align*}
V^i(\Gamma\times_{\Delta}\Lambda) &= \{(v, w)\in V^i(\Gamma)\times V^i(\Lambda) \mid \gamma(v) = \lambda(w)\},\\
V^f(\Gamma\times_{\Delta}\Lambda) &= \{(v, w)\in V^f(\Gamma)\times V^f(\Lambda) \mid \gamma(v) = \lambda(w)\}.
\end{align*}
In general, as pointed out in \cite{domaratzki_07}, the following is satisfied:
\[
\abs{E(\Gamma\times_{\Delta}\Lambda)} = \sum_{e\in E(\Delta)}\abs{\gamma_E^{-1}(e)}\cdot \abs{\lambda_E^{-1}(e)}.
\]
If $\gamma$ and $\lambda$ are finite covers, then $\Gamma\times_{\Delta}\Lambda$ is equal to its core. Hence, a quadratic bound on the edges and vertices is best possible. The following lemma demonstrates a similar relationship between $\abs{\bar{E}(\Gamma)}$, $\abs{\bar{E}(\Lambda)}$ and $\abs{\bar{E}(\rcore(\Gamma\times_{\Delta}\Lambda))}$ when $\gamma$ and $\lambda$ are forwards immersions.

\begin{lemma}
\label{ve_bounds}
Let $\gamma:\Gamma\to \Delta$ and $\lambda:\Lambda\to \Delta$ be forwards immersions of graphs. If $\Theta = \rcore(\Gamma\times_{\Delta}\Lambda)$, then
\[
\abs{\bar{V}(\Theta)}, \abs{\bar{E}(\Theta) - \bar{E}_c(\Theta)} \leq 4\cdot \left(\abs{V^i(\Gamma)}\cdot \abs{V^i(\Lambda)} + \abs{\bar{E}(\Gamma)}\cdot \abs{\bar{E}(\Lambda)}\right).
\]
\end{lemma}

\begin{proof}
Let $V'\subset V(\Theta)$ be the vertices with outdegree two or more. By removing these vertices and taking the closure, we obtain a disjoint union of graphs $\sqcup_{i=1}^n\Theta_i$ along with a surjective quotient map $q:\sqcup_{i=1}^n\Theta_i\to \Theta$. Each graph $\Theta_i$ has the property that $\odeg(v)\leq 1$ by construction. Furthermore, $q:\bar{V}(\sqcup_{i=1}^n\Theta_i)\to\bar{V}(\Theta)$ is surjective and $q:\bar{E}(\sqcup_{i=1}^n\Theta_i)\to \bar{E}(\Theta)$ is bijective. We have at most
\[
\abs{V^i(\Theta)} + \sum_{v\in V'}\odeg(v)\leq \abs{V^i(\Gamma)}\cdot \abs{V^i(\Lambda)} + \abs{\bar{E}(\Gamma)}\cdot \abs{\bar{E}(\Lambda)}
\]
many leaf vertices in $\sqcup_{i=1}^n\Theta_i$. Hence
\[
\sum_{i=1}^n \abs{\partial \Theta_i}\leq 2\cdot \left(\abs{V^i(\Gamma)}\cdot \abs{V^i(\Lambda)} + \abs{\bar{E}(\Gamma)}\cdot \abs{\bar{E}(\Lambda)}\right).
\]
Thus, since $\beta_1(\Theta_i)\leq 1$ for all $i$, the bounds follow by Lemma \ref{betti_bound}.
\end{proof}

\section{The structure of $\Gamma\times_{\Delta}S^1$}
\label{w_graphs_section}

Let us fix a graph map $\gamma:\Gamma\to \Delta$ and a primitive loop $w:I\to \Delta$. A cycle $f:S^1\to \Gamma$ is a \emph{$w$-cycle} if $f$ factors through $\hat{w}$. A path $f:I\to \Gamma$ will be called a \emph{$w$-path} if $\gamma\circ f$ factors surjectively through $\hat{w}$. The projection of the $w$-path to $\Delta$ can be factorised like a word:
\[
\gamma\circ f = w'*w^n*w''
\]
where $w'$ and $w''$ are a proper suffix and a proper prefix of $w$ respectively. We call this a \emph{$w$-factorisation of $f$}. The following Lemma follows from Lemma \ref{unique_factorisation}:

\begin{lemma}
\label{unique_w}
If $f:I\to \Gamma$ is a $w$-path, then it has a unique $w$-factorisation.

\end{lemma}

The \emph{$w$-length} of a $w$-path $f:I\to \Gamma$ is the number $n = |f|_w$ such that $\gamma\circ f = w'*w^n*w''$ is the $w$-factorisation of $f$. By Lemma \ref{unique_w} this is always well defined.

A $w$-path $f'$ \emph{right extends} $f$, or is a \emph{right extension of $f$}, if $f$ is a proper prefix of $f$. If $f$ is a proper suffix of $f'$, then we say $f'$ \emph{left extends} $f$. A $w$-path is a \emph{simple $w$-path} if 
\[
f(k + i\cdot |w|)\neq f(k + j\cdot |w|)
\]
for all integers $i< j$ and $k$, unless $k+j\cdot |w| = |f|$. We refer to $f$ as \emph{right maximal} if $f$ is a simple $w$-path and there is no simple $w$-path $f'$ which right extends $f$. We similarly define \emph{left maximal}. We refer to $f$ as \emph{maximal} if it is both left and right maximal. We also refer to $f$ as \emph{right submaximal} if either $f$ is right maximal, or $f$ is simple and $t(f)\in \bar{V}(\Gamma)\cup V^f(\Gamma)$. Similarly, we say $f$ is \emph{left submaximal} if either $f$ is left maximal, or if $f$ is a simple $w$-path and $o(f) \in \bar{V}(\Gamma)\cup V^i(\Gamma)$. Then \emph{submaximal} $w$-paths are those that are both left and right submaximal. Denote by $\wmax(\gamma, w)$ the set of maximal $w$-paths in $\Gamma$ and by $\subwmax(\gamma, w)$ the set of submaximal $w$-paths in $\Gamma$. 
Defining:
\begin{align*}
\swmax(\gamma, w) &= \slift(\gamma, w)\cap \wmax(\gamma, w), &\owmax(\gamma, w) &= \olift(\gamma, w)\cap \wmax(\gamma, w),\\
\ssubwmax(\gamma, w) &= \slift(\gamma, w)\cap \subwmax(\gamma, w), &\osubwmax(\gamma, w) &= \olift(\gamma, w)\cap \subwmax(\gamma, w),
\end{align*}
we get decompositions 
\begin{align*}
\wmax(\gamma, w) &= \swmax(\gamma, w)\sqcup\owmax(\gamma, w),\\
\subwmax(\gamma, w) &= \ssubwmax(\gamma, w)\sqcup \osubwmax(\gamma, w).
\end{align*}

\begin{remark}
If $\gamma:\Gamma\to \Delta$ is a forwards immersion, then maximal right extensions of $w$-paths are unique. However, there may be several left maximal extensions if $\gamma$ is not an immersion.
\end{remark}

The following lemmas follow directly from the definitions, but we record them for later use.

\begin{lemma}
\label{w-decomposition}
Let $\gamma:\Gamma\to \Delta$ be a forwards immersion and $w:I\to \Delta$ a primitive loop. Every $w$-path $f:I\to \Gamma$ has a unique decomposition:
\[
f = f_s*f_c^n*f_c'
\]
for some $n\geq 0$ and where
\begin{enumerate}
\item $f_s*f_c$ is a simple $w$-path, 
\item $f_c$ is a simple $w$-cycle,
\item $f_c'$ is a proper prefix of $f_c$.
\end{enumerate}
Furthermore, $f$ is simple if and only if $n\leq 1$ and $f_c' = \epsilon$.

\end{lemma}

\begin{lemma}
\label{submax_subpath}
Let $\gamma:\Gamma\to \Delta$ be a graph map and $w:I\to \Delta$ a primitive loop. Let $f:I\to \Gamma$ be an admissible path with a factorisation
\[
f = f_1*f_w*f_2
\]
with $f_w$ a maximal $w$-path in $I$. Then if $f_w = f_s*f_c^n*f_c'$ is the factorisation from Lemma \ref{w-decomposition}, $f_s*f_c$ is submaximal and $f_c'$ is a submaximal prefix of $f_c$.

\end{lemma}

Extending our definition of $w$-paths to graphs, we say a graph map $W\to \Delta$ is a \emph{$w$-graph map} if it factors surjectively through the map $\hat{w}:S^1\to \Delta$. The graph $W$ will then be called a \emph{$w$-graph}.

\begin{lemma}
\label{w-graphs_shape}
Let $W$ be a connected graph and $W\to S^1$ a forwards immersion. Then $W$ is one of the following:
\begin{enumerate}
\item A rooted tree,
\item a cycle with rooted trees attached along their root vertices.
\end{enumerate}
\end{lemma}

Every $w$-path factors through $\Gamma\times_{\Delta}S^1$, but this is not quite a $w$-graph as it contains components that don't surject onto $S^1$. So denote by $\Omega_{\gamma, w}\subset \Gamma\times_{\Delta}S^1$ the subgraph consisting of the $w$-graph components. Call this \emph{the $w$-graph of $\Gamma$} and denote by $p_w:\Omega_{\gamma, w}\to S^1$ the projection map. We will further define the subgraph $\Omega_{\gamma, w}^k\subset \Omega_{\gamma, w}$ consisting of the components which support a $w$-path of $w$-length at least $k$.

The following Lemma is a graph version of Lemma \ref{unique_factorisation}.

\begin{lemma}
\label{single_lift}
Let $f:I\to \Gamma$ be a path with $|f|\geq |w|$ and $p_{\Gamma}:\Omega_{\gamma, w}\to \Gamma$ the projection map, then:
\[
\abs{\lift(p_{\Gamma}, f)}\leq 1.
\]
\end{lemma}

\begin{proof}
If $\tilde{f}:I\to\Omega_{\gamma, w}$ is a lift of $f$, then $f$ has a $w$-factorisation. If $|f|\geq |w|$, by Lemma \ref{unique_factorisation}, $f$ has a unique $w$-factorisation and so $o(p_{w}\circ \tilde{f})$ is uniquely determined by this factorisation. By the construction of the fibre product graph, it follows that there can be at most one lift of $f$ to $\Omega_{\gamma, w}$.
\end{proof}

\subsection{$w$-paths}

In this section we will fix a graph map $\gamma:\Gamma\to \Delta$ and a primitive loop $w:I\to \Delta$. The aim of this section will be to prove bounds on the number of $w$-lengths that $w$-paths in $\Gamma$ can have.

First, we define a type of factorisation of paths relative to a subgraph.

\begin{definition}
Let $\Gamma'\subset \Gamma$ be a subgraph and $f:I\to \Gamma$ a path. If $f$ does not factor through $\Gamma - \bar{V}(\Gamma)$, the \emph{$\Gamma'$-factorisation} of $f$ is the factorisation of the form:
\[
f = e_0[-l:]*w_1*e_1*...*e_{n-1}*w_n*e_n[:m],
\]
where:
\begin{enumerate}
\item $o(w_i)$, $t(w_i)\in \bar{V}(\Gamma)$,
\item $w_i$ is supported in $\Gamma'\cup \bar{V}(\Gamma)$,
\item $e_0$, $e_n\in \bar{E}(\Gamma)$,
\item $e_1, ..., e_{n-1}\in \bar{E}(\Gamma) - \bar{E}(\Gamma')$,
\item $n$ is maximal possible among all such factorisations.
\end{enumerate}
If $f$ does factor through $\Gamma - \bar{V}(\Gamma)$, then by convention the $\Gamma'$-factorisation of $f$ will be $e_0[l:m]$ for some $e_0\in \bar{E}(\Gamma)$.
\end{definition}

The following lemma follows from the definition.

\begin{lemma}
$\Gamma'$-factorisations are unique.

\end{lemma}

Denote by $\bar{E}_w(\Gamma)\subset \bar{E}(\Gamma)$ the elements $e$ such that $|e|<|w|$. Then denote by
\[
\Gamma_w = \cup_{e\in\bar{E}_w(\Gamma)}e
\]
the subgraph of $\Gamma$ induced by $\bar{E}_w(\Gamma)$. 

\begin{definition}
We will call a vertex $v\in V(\Gamma)$ a \emph{$w$-sink} if there is a $w$-path $f:I\to \Gamma_w$ satisfying the following:
\begin{enumerate}
\item $o(f)$, $t(f) \in \bar{V}(\Gamma)$,
\item $\gamma\circ f = w'*w*w''$ for some proper suffix $w'$ of $w$ and proper prefix $w''$ of $w$,
\item $f(|w'*w|) = v$.
\end{enumerate}
Denote by $V_w(\Gamma)\subset V(\Gamma_w)$ the set of all $w$-sinks of $\Gamma$.
\end{definition}

The $w$-sinks of a graph can be thought of as a set of vertices any sufficiently long $w$-path must visit. We now show that there cannot be too many of such vertices.

\begin{lemma}
\label{general_w_sinks}
There is a cyclic shift $w_0$ of $w$ such that:
\[
\abs{V_{w_0}(\Gamma)}\leq 2\cdot \abs{\bar{E}(\Gamma_w)}.
\]
\end{lemma}

\begin{proof}
Since a cycle doesn't have any $w$-sinks, we may assume that $\bar{E}_c(\Gamma) = \emptyset$. Now let $w = x*y$ with $x\neq \epsilon$ such that $w_0 = y*x$ is the cyclic shift of $w$ from Theorem \ref{cyclic_shift} with 
\[
\{x_1, ..., x_n\} = \{\gamma\circ e \mid e\in \bar{E}(\Gamma_w)\}.
\]
For each $v\in V_{w_0}(\Gamma)$, let $f_v:I\to \Gamma_w$ be a $w_0$-path such that: $o(f_v)$, $t(f_v)\in \bar{V}(\Gamma_w)$, $\gamma\circ f_v = w_v'*w_0*w_v''$ and $f_v(|w_v'*w_0|) = v$. Denote by
\[
\bar{t}\left(f_v\left[:\abs{w_v'*w_0}\right]\right) = (e_v, i_v)\in \bar{E}(\Gamma_w)\times \N
\]
Now by Theorem \ref{cyclic_shift}, we have:
\[
\abs{\{(e_v, i_v)\}_{v\in V_{w_0}(\Gamma)}\cap \{e\}\times\N}\leq 2,
\]
for all $e\in \bar{E}(\Gamma_w)$, thus obtaining the desired inequality.
\end{proof}

\begin{lemma}
\label{w-sinks}
If $\gamma:\Gamma\to \Delta$ is a forwards immersion, then:
\[
\abs{V_w(\Gamma)}\leq 2\cdot \abs{\bar{E}(\Gamma_w)}.
\]
\end{lemma}

\begin{proof}
Since a cycle doesn't have any $w$-sinks, we may assume that $\bar{E}_c(\Gamma) = \emptyset$. Now let $w = x*y$ with $x\neq \epsilon$ such that $y*x$ is the cyclic shift of $w$ from Theorem \ref{cyclic_shift} with 
\[
\{x_1, ..., x_n\} = \{\gamma\circ e \mid e\in \bar{E}(\Gamma_w)\}.
\]
For each $v\in V_w$, let $f_v:I\to \Gamma_w$ be a $w$-path such that: $o(f_v)$, $t(f_v)\in \bar{V}(\Gamma)$, $\gamma\circ f_v = w_v'*w*w_v''$ and $f_v(|w_v'*w|) = v$. Denote by
\[
\bar{t}(f_v[:\abs{w_v'*x}]) = (e_v, i_v)\in \bar{E}(\Gamma_w)\times \N
\]
This is well defined since $x\neq \epsilon$. Since $\gamma$ is a forwards immersion, $(e_v, i_v) = (e_{v'}, i_{v'})$ if and only if $v = v'$. Now by Theorem \ref{cyclic_shift}, we have:
\[
\abs{\{(e_v, i_v)\}_{v\in V_w(\Gamma)}\cap \{e\}\times\N}\leq 2,
\]
for all $e\in \bar{E}(\Gamma_w)$, thus obtaining the desired inequality.
\end{proof}

Using these bounds on $w$-sinks and $\Gamma_w$-factorisations, we prove something about general $w$-paths.

\begin{lemma}
\label{gamma_w_factorisation}
Let $f:I\to \Gamma$ be a simple $w$-path that does not factor through $\Gamma - \bar{V}(\Gamma)$. Let
\[
f = e_0[-l:]*w_1*e_1*...*e_{n-1}*w_n*e_n[:m],
\]
be its $\Gamma_w$-factorisation. Then:
\begin{enumerate}
\item\label{itm:equal} $e_i = e_j$ if and only if $i = j$ for all $1\leq i, j<n$,
\item\label{itm:nequal} if $l\geq \abs{w}$ or $m\geq \abs{w}$, then $e_0\neq e_i$ or $e_m\neq e_i$ for all $1\leq i<n$ respectively,
\item\label{itm:sum} $\sum_{i=1}^n \abs{w_i}_w\leq 4\cdot \abs{\bar{E}(\Gamma_w)}$,
\item\label{itm:fsum} if $\gamma$ is a forwards immersion, then $\sum_{i=1}^n \abs{w_i}_w\leq 2\cdot \abs{\bar{E}(\Gamma_w)}$
\end{enumerate}
\end{lemma}

\begin{proof}
Since $f$ is a simple $w$-path, it lifts to $\Omega_{\gamma, w}$ as a simple path. As the lift is simple as a path, it can only traverse each edge in $\Omega_{\gamma, w}$ at most once. Hence by Lemma \ref{single_lift}, \ref{itm:equal} and \ref{itm:nequal} follow. Property \ref{itm:sum} follows from Lemma \ref{general_w_sinks}, the pigeonhole principle and the fact that $w^2$ contains a copy of each cyclic shift of $w$ as a subword. More directly, \ref{itm:fsum} follows from Lemma \ref{w-sinks}.
\end{proof}

A direct consequence of Lemma \ref{gamma_w_factorisation} is that if $f:I\to \Gamma$ is a simple $w$-path, then:
\[
\abs{f}_w\leq 4\cdot \abs{\bar{E}(\Gamma)} + \sum_{e\in \bar{E}(\Gamma)}\left\lfloor\frac{\abs{e}}{\abs{w}}\right\rfloor.
\]
The following theorem is an important ingredient in our main results.

\begin{theorem}
\label{right_maximal_integers}
Suppose $\gamma:\Gamma\to \Delta$ is a forwards immersion. There are at most $10\cdot \abs{\bar{E}(\Gamma)} + 3\cdot \abs{V^i(\Gamma)}$ integers $n\geq 1$ such that there is a right maximal $w$-path $f\in \osubwmax(\gamma, w)$ with $|f|_w = n$.
\end{theorem}

\begin{proof}
We may assume that $\Gamma$ is connected. If $\abs{\bar{E}(\Gamma)} = 1$, then Lemma \ref{unique_factorisation} tells us there can be at most $2 + \abs{V^i(\Gamma)}$ such integers $n$: one integer for the maximal $w$-prefix of $e\in \bar{E}(\Gamma)$, one integer for the maximal $w$-suffix of $e$ and one integer for each vertex in $V^i(\Gamma)$ from which a $w$-path may begin. So from now on we assume that $\abs{\bar{E}(\Gamma)}\geq 2$. Let $f:I\to \Gamma$ be a right maximal submaximal $w$-path and let
\[
f = e_0[-l:]*w_1*e_1*...*w_n*e_n[:m]
\]
be its $\Gamma_w$-factorisation. If $f$ is a loop, we may assume that $o(f), t(f)\in \bar{V}(\Gamma)$. Since $f$ is a $w$-path, each subpath $e_i$ and $w_i$ has a well defined $w$-factorisation. So let:
\begin{align*}
\gamma\circ e_0[-l:] &= w_{e_0}'*w^{q_0}*w_{e_0}''\\
\gamma\circ e_i &= w_{e_i}'*w^{q_i}*w_{e_i}''\\
\gamma\circ e_n[:m] &= w_{e_n}'*w^{q_n}*w_{e_n}''\\
\gamma\circ w_i &= w_i'*w^{k_i}*w_i'' 
\end{align*}
be the relevant $w$-factorisations. Denote by $V_f\subset V(\Gamma)$ the set of vertices that fall into one of the following categories:
\begin{enumerate}
\item\label{itm:e_i} $v = e_i(\abs{w_{e_i}'})$ for some $1\leq i< n$, or,
\item\label{itm:e_0} $v = e_0(-l + \abs{w_0'})$ if $l\geq |w|$, or, 
\item\label{itm:e_n} $v = e_n(\abs{w_{e_n}'})$ if $q_n\geq 1$, or,
\item\label{itm:w_i} $v = w_i(\abs{w_i'*w^j})$ for some $i\geq 1$ and $1\leq j\leq k_i$, or,
\item\label{itm:w_1} $v = w_1(\abs{w_1'})$ if $l \geq |w|$, or,
\item\label{itm:w_2} $v = w_i(\abs{w_i'})$ for some $i\geq 2$.
\end{enumerate}
By the definition of $\Gamma_w$-factorisations, each vertex of type (\ref{itm:e_i}) and (\ref{itm:e_n}) is contained in some element $e\in \bar{E}(\Gamma) - \bar{E}_w(\Gamma)$. Moreover, these are independent of $f$ by Lemma \ref{unique_factorisation}. Hence there can be at most $\abs{\bar{E}(\Gamma) - \bar{E}_w(\Gamma)}$ many vertices $v\in \bigcup_f V_f$ of type (\ref{itm:e_i}) or (\ref{itm:e_n}), where the union is over all right maximal submaximal $w$-paths. Since either $o(e_0[-l:])\in V^i(\Gamma)$, or $e_0[-l:]$ is a left maximal $w$-path, there can be at most $\abs{\bar{E}(\Gamma) - \bar{E}_w(\Gamma)} + \abs{V^i(\Gamma)}$ many vertices of type (\ref{itm:e_0}). By Lemma \ref{w-sinks}, there are at most $2\cdot \abs{\bar{E}_w(\Gamma)}$ vertices of type (\ref{itm:w_i}) in $\bigcup_f V_f$. Since right extensions of $w$-paths are uniquely defined for forwards immersions, we get a bound of $\abs{\bar{E}(\Gamma) - \bar{E}_w(\Gamma)}$ for the vertices of type (\ref{itm:w_1}) and (\ref{itm:w_2}).

So putting all of this together we get:
\[
\abs{\bigcup_fV_f}\leq 3\cdot \abs{\bar{E}(\Gamma)} + \abs{V^i(\Gamma)}.
\]
To each vertex $v\in \bigcup_f V_f$, we may associate an integer $k_v$ that is the $w$-length of the right maximal $w$-path starting at $v$ with $w$-factorisation beginning with $w$.

Now let $g:I\to \Gamma$ be a right maximal submaximal $w$-path and let $g = w'*w^k*w''$ be its $w$-factorisation. Then by definition of $V_g\subset \bigcup_f V_f$, either $g(\abs{w'})\in V_g$, $g(\abs{w'*w})\in V_g$, $g(\abs{w'*w^2})\in V_g$ or $\abs{g}_w\leq 2$. Hence $|g|_w = k_v$, $|g|_w = k_v + 1$ or $|g|_w = k_v + 2$ for some $v\in \bigcup_f V_f$, or $|g|_w\leq 2$. So $k$ can be one of at most $3\cdot \left(3\cdot \abs{\bar{E}(\Gamma)} + \abs{V^i(\Gamma)}\right) + 2$ distinct integers. Since we assumed that $\abs{\bar{E}(\Gamma)}\geq 2$, we obtain the desired upper bound.
\end{proof}

\begin{remark}
If we wanted to relax the assumptions of Theorem \ref{right_maximal_integers} to encompass general graph maps, then the number of integers $n$ needed could be exponential in $\abs{\bar{E}(\Gamma)}$.
\end{remark}

\begin{theorem}
\label{extension_language}
Suppose that $\gamma:\Gamma\to \Delta$ is a forwards immersion and let $(u, v)\in V(\Gamma)\times V(S^1)$. There is a regular language of the form:
\[
L = (a^{n_1} + ... + a^{n_k})\cdot (a^n)^*\subset (a)^*
\]
such that, for any path $f:I\to \Omega_{\gamma, w}$ where $t(f) = (u, v)$ and $p_{\Gamma}\circ f$ is a left submaximal $w$-path, we have $a^{|f|_w}\in L$. Furthermore, we have that:
\[
k \leq 10\cdot \abs{\bar{E}(\Gamma)} + 3\cdot \abs{V^i(\Gamma)}.
\]
\end{theorem}

\begin{proof}
The fact that the regular language is of the given form is by Lemma \ref{w-graphs_shape}. If there is an $f:I\to \Omega_{\gamma, w}$ such that $t(f) = (u, v)$, $p_{\Gamma}\circ f$ is a left submaximal $w$-path and $p_{\Gamma}\circ f$ does not contain any vertex in $\bar{V}(\Gamma)$ in its image, then $k = 1$ and $n= 0$. So, assuming that no such path exists, the proof of the bound on the number of integers is identical to that of Theorem \ref{right_maximal_integers}.
\end{proof}

\begin{remark}
We have stated Theorems \ref{right_maximal_integers} and \ref{extension_language} for left submaximal $w$-paths, but the analogous results for right submaximal paths also hold with virtually the same proof if we replace $V^i(\Gamma)$ with $V^f(\Gamma)$.
\end{remark}

\subsection{$w$-cycles}

The following theorem was conjectured in \cite{wise_05} in connection with coherence of one-relator groups and became known as the $w$-cycles conjecture. It was then proved independently in \cite{helfer_16} and \cite{louder_17}.

\begin{theorem}
\label{w_cycles}
Let $\gamma:\Gamma\to \Delta$ be an immersion and $\hat{w}:S^1\to \Delta$ a primitive cycle. Then:
\[
\beta_0(\core(\Gamma\times_{\Delta}S^1))\leq \beta_1(\Gamma).
\]
\end{theorem}

If we want to allow $\hat{w}$ to be imprimitive, then it follows that the map $\core(\Gamma\times_{\Delta}S^1)\to \Delta$ consists of at most $\beta_1(\Gamma)$ many free homotopy classes of cycles. We here prove a result transverse to this, Theorem \ref{long_cycles_graph_equiv} from the introduction.

\begin{lemma}
\label{long_w_edge}
Let $\gamma:\Gamma\to \Delta$ be a graph map and $w:I\to \Delta$ a primitive loop. Let $f:S^1\to \Gamma$ be a primitive cycle with $\deg(\gamma\circ f)\geq 12\cdot \beta_1(\Gamma)$. Then there is some $e\in \bar{E}(\Gamma)$ such that $f$ traverses $e$ and $\abs{e}\geq 2\cdot \abs{w}$.
\end{lemma}

\begin{proof}
We may assume that $\Gamma$ is connected and core. If $\Gamma$ is a cycle then the result is clear, so suppose that it is not. By Lemma \ref{betti_bound}, $\abs{\bar{E}(\Gamma)}\leq 3\cdot (\beta_1(\Gamma) - 1)$. Let $g:I\to \Gamma$, $w:I\to \Delta$ be loops such that $\hat{g} = f$ and $\gamma\circ g = w^{\deg(f)}$. Then $\abs{g}_w>4\cdot \abs{\bar{E}(\Gamma)}$. Now we simply apply the pigeonhole principle to Lemma \ref{gamma_w_factorisation}.
\end{proof}

\begin{theorem}
\label{long_cycles_graph}
Let $\gamma:\Gamma\to \Delta$ be a forwards immersion of graphs. There is a unique graph map $s_{\Gamma}:\mathbb{S}_{\Gamma}\to \Gamma$ where $\mathbb{S}_{\Gamma}$ is a disjoint union of at most $\beta_1(\Gamma)$ circles such that, if $f:S^1\to \Gamma$ is a cycle satisfying:
\[
\frac{\deg(\gamma\circ f)}{\deg(f)}\geq 12\cdot \beta_1(\Gamma),
\]
then $f$ factors through $s_{\Gamma}$.
\end{theorem}

\begin{proof}
Each such cycle factors through a primitive cycle $g:S^1\to \Gamma$ such that $\deg(\gamma\circ g)\geq 12\cdot \beta_1(\Gamma)$. Hence, let $s_{\Gamma}:\mathbb{S}_{\Gamma}\to \Gamma$ be the union over all such cycles. Let $f_1, ..., f_n:S^1\to \Gamma$ be the cycles in $\mathbb{S}_{\Gamma}$. Let $w_i:I\to \Delta$ be a primitive loop such that $f_i$ is a $w_i$-cycle. Assume that $\abs{w_i}\geq \abs{w_{i+1}}$. By Lemma \ref{long_w_edge}, for each $f_i$ there is an element $e_i\in \bar{E}(\Gamma)$ such that $e_i$ is a $w_i$-path with $\abs{e_i}\geq 2\cdot \abs{w_i}$. If $f_{j}$ traverses $e_i$ for any $i<j$, up to replacing $w_i$ with a cyclic conjugate of itself, $w_i^2$ has period $\abs{w_j}\leq \abs{w_i}$. Since $w_i$ and $w_j$ are primitive, it follows from Corollary \ref{square_period} that $w_i$ is a cyclic conjugate of $w_j$. But since $\gamma$ was a forwards immersion, it then also follows that $f_i = f_j$. Hence $f_i$ is supported in $\Gamma - \{e_1, ..., e_{i-1}\}$. As $\Gamma - \{e_1, ..., e_{\beta_1(\Gamma)}\}$ is a forest, we must have $n\leq \beta_1(\Gamma)$. 
\end{proof}

Equivalently, we may restate Theorem \ref{long_cycles_graph} as follows.

\begin{theorem}
\label{long_cycles_graph_equiv}
Let $\gamma:\Gamma \to \Delta$ and $\lambda:\mathbb{S}\to \Delta$ be graph immersions with $\mathbb{S} = \sqcup S^1$. If $\lambda|_{S^1}$ is primitive for each $S^1\subset \mathbb{S}$, then there are at most $\beta_1(\Gamma)$ many components $S^1\subset \mathbb{S}$ such that:
\[
\deg(p_{\mathbb{S}}|_{S^1})\geq 12\cdot \beta_1(\Gamma).
\]

\end{theorem}

\subsection{$w$-graphs}

The following results can also be thought of as generalisations of Theorem \ref{w_cycles} as we do not only consider the core of the fibre product. Furthermore, the graph maps no longer need to be immersions.

\begin{theorem}
\label{w_components}
Let $\gamma:\Gamma\to \Delta$ be a forwards immersion and $\hat{w}:S^1\to \Delta$ a primitive cycle. There are at most $2\cdot \abs{\bar{E}(\Gamma)}$ components of $\Omega_{\gamma, w}^3$ that do not factor through $\Gamma - \bar{V}(\Gamma)$.
\end{theorem}

\begin{proof}
Let $\{f_i:I\to \Gamma\}_{i=1}^k$ be a collection of $w$-paths containing at most one of maximal $w$-length for each component of $\Omega_{\gamma, w}^3$ that does not factor through $\Gamma - \bar{V}$. Consider the $\Gamma_w$-factorisations for each 
\[
f_i = e_{i, 0}[-l_i:]*w_{i, 1}*e_{i, 1}...*w_{i, n_i}*e_{i, n_i}[:m_i].
\]
Since $f_i$ does not factor through $\Gamma - \bar{V}(\Gamma)$, $n_i\geq 1$ for all $i$. Since $\abs{f_i}_w\geq 3$ for all $i$, we may divide our set of paths into three types:
\begin{enumerate}
\item\label{itm:greater_one} those with $n_i>1$,
\item\label{itm:equal_one} those with $n_i = 1$ and $|w_{i, 1}|_w\geq 1$,
\item\label{itm:equal_one_2} those with $n_i = 1$, $|w_{i, 1}|_w = 0$ and $|e_{i, 0}[-l_i:]|_w\geq 1$ or $|e_{i, 1}[:m_i]|_w \geq 1$.
\end{enumerate}
By Lemma \ref{single_lift}, there can be at most $2\cdot \abs{\bar{E}(\Gamma) - \bar{E}_w(\Gamma)}$ paths $f_i$ in case \ref{itm:greater_one} or \ref{itm:equal_one_2}. By Lemma \ref{w-sinks} there can be at most $2\cdot \abs{\bar{E}_w(\Gamma)}$ many paths $f_i$ in case \ref{itm:equal_one}.
\end{proof}

\begin{theorem}
\label{w_components_general}
Let $\gamma:\Gamma\to \Delta$ be a graph map and $\hat{w}:S^1\to \Delta$ a primitive cycle. There are at most $2\cdot \abs{\bar{E}(\Gamma)}$ components of $\Omega_{\gamma, w}^4$ that don't factor through $\Gamma - \bar{V}(\Gamma)$.
\end{theorem}

\begin{proof}
Let $w = x*y$ with $x\neq \epsilon$ such that $w_0 = y*x$ is the cyclic shift of $w$ from Lemma \ref{general_w_sinks}. Now if we replace the use of Lemma \ref{w-sinks} in the proof of Theorem \ref{w_components} for Lemma \ref{general_w_sinks}, we get the same bound for $\Omega^3_{\gamma, w_0}$. Since each component of $\Omega^4_{\gamma, w}$ is a component of $\Omega^3_{\gamma, w_0}$, we obtain the bound.
\end{proof}

Theorems \ref{w_components} and \ref{w_components_general} are sharp as is demonstrated by Example \ref{sharp_components}. Obtaining similar statements for $\Omega_{\gamma, w}$ appears harder.

\begin{figure}
\centering
\vspace*{8pt}
\includegraphics[scale=0.7]{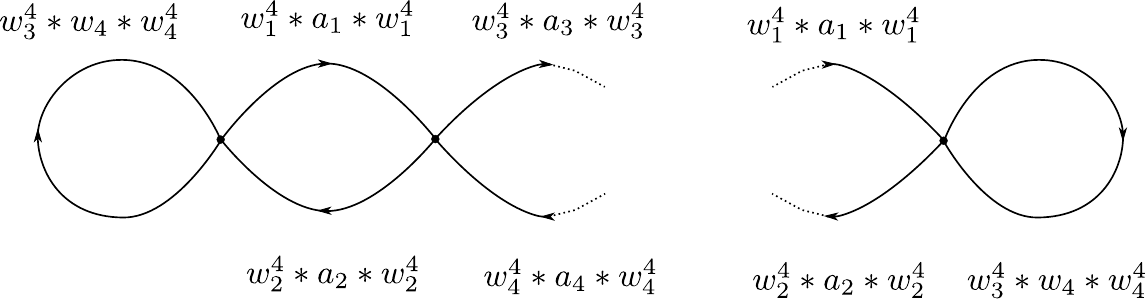}
\vspace{8pt}
\caption{When $n$ is even.}
\label{sharp_example}
\end{figure}

\begin{figure}
\centering
\vspace*{8pt}
\includegraphics[scale=0.7]{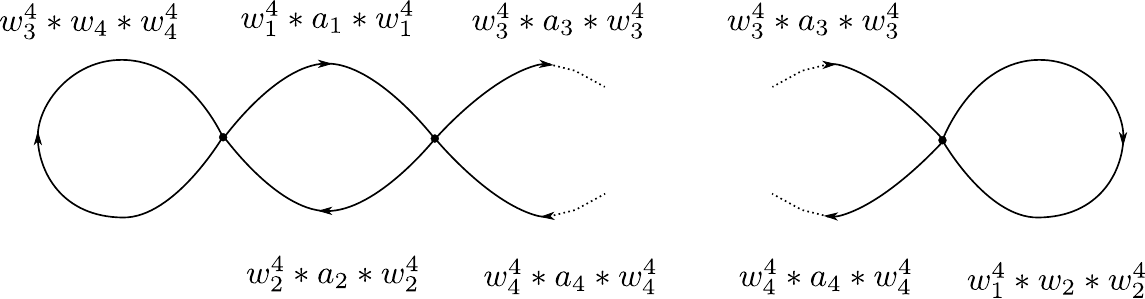}
\vspace{8pt}
\caption{When $n$ is odd.}
\label{sharp_example_2}
\end{figure}

\begin{example}
\label{sharp_components}
Let $\Delta$ be the graph with four edges and one vertex. Let $E(\Delta) = \{a_1, a_2, a_3, a_4\}$ and $w_i = a_i*a_{i+1}*a_{i+2}*a_{i+3}$ where the indices are taken mod $4$. Now for each integer $n\geq 1$, consider the graph immersion $\gamma:\Gamma\to\Delta$ with $\abs{\bar{V}(\Gamma)} = n$ depicted in Figure \ref{sharp_example} when $n$ is even and in Figure \ref{sharp_example_2} when $n$ is odd. We have $\abs{\bar{E}(\Gamma)} = 2\cdot n$ and $\Omega_{\gamma, w_1}^3$ has precisely $4\cdot n$ many components, none of which factor through $\Gamma - \bar{V}(\Gamma)$.
\end{example}

\section{The structure of $\Gamma\times_{\Delta}\Lambda$}
\label{main_structure_section}

From now on we will always assume that $\gamma:\Gamma\to \Delta$ and $\lambda:\Lambda\to \Delta$ are forwards immersions of graphs. For convenience we also assume that $\Gamma$ and $\Lambda$ are connected and relatively core. Let $s_{\Gamma}:\mathbb{S}_{\Gamma}\to \Gamma$ and $s_{\Lambda}:\mathbb{S}_{\Lambda}\to\Lambda$ be the graphs from Theorem \ref{long_cycles_graph}. We will refer to $\mathbb{S}_{\Gamma}$ and $\mathbb{S}_{\Lambda}$ as \emph{the long cycles} of $\Gamma$ and $\Lambda$ respectively. Then we define:
\begin{align*}
\theta:\Theta &= \rcore(\Gamma\times_{\Delta}\Lambda)\to \Delta,\\
s:\mathbb{S} &= \core(\mathbb{S}_{\Gamma}\times_{\Delta}\mathbb{S}_{\Lambda})\to \Theta.
\end{align*}
Define $\Theta^s\subset \Theta$ to be the subgraph with edge set 
\[
E(\Theta^s) = E(\Theta) - E(s(\mathbb{S}))
\]
and by
\begin{align*}
p_{\Gamma}^s&:\Theta^s\to \Gamma\\
p_{\Lambda}^s&:\Theta^s\to \Lambda.
\end{align*}
the restriction maps. Finally, denote by
\begin{align*}
C_{\Gamma} &= \abs{\bar{E}(\Gamma)} + \abs{V^i(\Gamma)},\\
C_{\Lambda} &= \abs{\bar{E}(\Lambda)} + \abs{V^i(\Lambda)}.
\end{align*}

\subsection{Counting lifts}
\label{counting_lifts}

We here apply our $w$-paths results to count certain subsets of lifts to the fibre product graph. The idea is that we will apply this to elements of $\bar{E}(\Gamma)$ and $\bar{E}(\Lambda)$ to cover the graph $\rcore(\Gamma\times_{\Delta}\Lambda)$. The reader should have in mind the connection between the following Proposition and the results in Section \ref{word_overlaps} on word overlaps.

\begin{proposition}
\label{same_period}
Let $g:I\to \Gamma$ be a simple path, $e\in \bar{E}(\Lambda) - \bar{E}_c(\Lambda)$ and $\ap(p, q, r)$ an arithmetic progression such that $p\cdot q\leq \abs{g}$. Then:
\[
\abs{\{\tilde{g} \in \lift(p^s_{\Gamma}, g) \mid \bar{o}(p^s_{\Lambda}\circ \tilde{g})\in \{e\}\times\ap(p, q, r)\}}\leq 100\cdot C_{\Gamma}\cdot C_{\Lambda}.
\]
\end{proposition}

\begin{proof}
Let $\{\tilde{g}_1, ..., \tilde{g}_n\} = \{\tilde{g}\in \lift(p^s_{\Gamma}, g) \mid \bar{o}(p_{\Lambda}\circ \tilde{g})\in \{e\}\times \ap(p, q, r)\}$ and let $p_1, ..., p_n$ be integers such that $\bar{o}(p_{\Lambda}\circ \tilde{g}_i) = (e, p_i)$. Note that since $\gamma$ and $\lambda$ are forwards immersions, each lift of $g$ is uniquely determined by its origin. Thus, $p_i\neq p_j$ for all $i\neq j$. We assume that the $\tilde{g}_i$ are ordered such that $p_1<...<p_n$. By assumption $p_n - p_1\leq \abs{g}$ and $q\mid p_j - p_i$ for all $i<j$. By Lemma \ref{key_power}, $e[p_1:p_n]$ has a period of $p_j - p_i$ for all $i<j$. If $n\leq 2$, then we are done so we may assume that $n\geq 3$. It follows that $e[p_1:p_n] = w^l$ for some $l\geq 2$ and $w$ primitive such that $\abs{w}\mid p_j - p_i$ for all $i<j$. Thus, we may assume without loss that $p = k$, $q = \abs{w}$ and $r = p_1$.

Let $e(p_i) = v_i$. For each $i$, choose a left maximal $w$-path $h_i:I\to \Theta$ such that $t(h_i) = (o(g), v_i)$. Note that $p_{\Lambda}\circ h_i$ and $p_{\Gamma}\circ h_i$ are both left submaximal $w$-paths by Lemma \ref{submax_subpath} for all $i\geq 1$. Let $m_1, ..., m_{\alpha}, m$ be the integers from Theorem \ref{extension_language} applied to $(o(g), o)\in V(\Gamma)\times V(S^1)$, where $o\in V(S^1)$ is the vertex from which we read off $w$. Similarly, let $k_1, ..., k_{\beta}, k$ be the integers from Theorem \ref{extension_language} when applied to $(e(r), o)\in V(\Lambda)\times V(S^1)$. We can see that for all $i\geq 1$, we must have:
\[
a^{\abs{h_i}_w} \in  ((a^{k_1} + ... + a^{k_{\beta}})\cdot (a^k)^*\cdot (\emptyset + a + ... + a^{l})) \cap ((a^{m_1} + ... + a^{m_{\alpha}})\cdot (a^m)^*).
\]
Since $g$ was a simple path, if $m\neq 0$, we must have $l\leq m$. If $k\neq 0$, then also $l\leq k$. By hypothesis, $\tilde{g}_i$ is supported in $\Theta^s$ for each $i$. Hence either $k < 12\cdot \beta_1(\Lambda)$ or $m < 12\cdot\beta_1(\Gamma)$ by construction. So if $k, m\neq 0$, then we are done. So now assume that $k= 0$ or $m=0$. This regular language can consist of at most $\alpha\cdot \beta$ elements. In particular, there are at most $\alpha\cdot \beta$ many integers $j$ such that 
\[
((a^{k_1} + ... + a^{k_{\beta}})\cdot(a^{k})^*\cdot a^j)\cap ((a^{m_1} + ... + a^{m_{\alpha}})\cdot (a^m)^*)\neq \emptyset.
\]
Thus, we have $n\leq \alpha\cdot \beta$. Finally, by Theorem \ref{extension_language}:
\[
n\leq \left(10\cdot \abs{\bar{E}(\Gamma)} + 3\cdot \abs{V^i(\Gamma)}\right)\cdot \left(10\cdot \abs{\bar{E}(\Lambda)} + 3\cdot \abs{V^i(\Lambda)}\right)\leq 100\cdot C_{\Gamma}\cdot C_{\Lambda}.
\]
\end{proof}

\subsection{A decomposition of $E(\Theta)$}
\label{decomposition}

Let $e\in \bar{E}(\Gamma)$ and $f\in \bar{E}(\Lambda)$. Denote by:
\[
\Theta(e, f) = p_{\Gamma}^{-1}(e)\cap p_{\Lambda}^{-1}(f)\subset \Theta.
\]
Then we have a decomposition of $E(\Theta)$ into disjoint sets of edges coming from each subgraph $\Theta(e, f)$:
\[
E(\Theta) = \bigsqcup_{e\in \bar{E}(\Gamma)}\bigsqcup_{f\in \bar{E}(\Lambda)}E(\Theta(e, f)).
\]

\begin{lemma}
\label{segment_component}
Each component of $\Theta(e, f)$ is a segment.
\end{lemma}

\begin{proof}
Each component of $\Theta(e, f)$ projects into $e$ and $f$. Hence, apart from the boundary vertices, each vertex must have indegree and outdegree precisely one.
\end{proof}

By Lemma \ref{segment_component} we may consider each component $g\subset \Theta(e, f)$ as a directed path $g:I\to \Theta$. Just as in Section \ref{word_overlaps}, we define subsets $\Theta_{\sub}(e, f), \Theta_{\super}(e, f), \Theta_{\pref}(e, f), \Theta_{\suff}(e, f), \Theta_{\partial}(e, f)\subset \Theta(e, f)$ in the following way. Let $g\subset \Theta(e, g)$ be a component, then:
\begin{align*}
g&\subset \Theta_{\sub}(e, f) && \text{if } o(p_{\Lambda}\circ g) = o(f) \text{ and } t(p_{\Lambda}\circ g) = t(f)\\
g&\subset \Theta_{\super}(e, f) && \text{if }o(p_{\Gamma}\circ g) = o(e) \text{ and } t(p_{\Gamma}\circ g) = t(e)\\
g&\subset \Theta_{\pref}(e, f)&& \text{if } o(p_{\Gamma}\circ g) = o(e) \text{ and } t(p_{\Lambda}\circ g) = t(f)\\
g&\subset \Theta_{\suff}(e, f) && \text{if } o(p_{\Lambda}\circ g) = o(f) \text{ and } t(p_{\Gamma}\circ g) = t(e)\\
g&\subset \Theta_{\partial}(e, f) && \text{if } o(g)\in \partial\Theta \text{ or } t(g)\in \partial\Theta
\end{align*}
These subsets provide us with a further decomposition:
\[
E(\Theta(e, f)) = E(\Theta_{\sub}(e, f))\cup E(\Theta_{\super}(e, f))\cup E(\Theta_{\pref}(e, f))\cup E(\Theta_{\suff}(e, f))\cup E(\Theta_{\partial}(e, f)).
\]
As we will treat the subgraph $\Theta^s$ separately, we will mostly be interested in the following subgraphs:
\begin{align*}
\Theta^s(e, f) &= \Theta(e, f)\cap \Theta^s,\\
\Theta_{\sub}^s(e, f) &= \Theta_{\sub}(e, f)\cap \Theta^s,\\
\Theta_{\super}^s(e, f) &= \Theta_{\super}(e, f)\cap \Theta^s,\\
\Theta_{\pref}^s(e, f) &= \Theta_{\pref}(e, f)\cap \Theta^s,\\
\Theta_{\suff}^s(e, f) &= \Theta_{\suff}(e, f)\cap \Theta^s.
\end{align*}
Thus, we get a decomposition of $E(\Theta)$.

\begin{lemma}
\label{edge_decomposition}
\begin{align*}
E(\Theta) = E(s(\mathbb{S})) \sqcup \bigsqcup_{e\in \bar{E}(\Gamma)}\bigsqcup_{f\in \bar{E}(\Lambda)}&(E(\Theta^s_{\sub}(e, f))\cup E(\Theta^s_{\super}(e, f)) \cup\\
		&\cup E(\Theta^s_{\pref}(e, f))\cup E(\Theta^s_{\suff}(e, f))\cup E(\Theta_{\partial}(e, f))).
\end{align*}

\end{lemma}

\subsection{Bounding $|E(\Theta^s)|$}

In this section we will show that each of the sets $E(\Theta_{\sub}^s(e, f))$, $E(\Theta^s_{\super}(e, f))$, $E(\Theta^s_{\pref}(e, f))$, $E(\Theta^s_{\suff}(e, f))$ and $E(\Theta_{\partial}(e, f))$ have, up to multiplication by some universal constant, cardinality bounded above by $C_{\Gamma}\cdot C_{\Lambda}\cdot (\abs{e} + \abs{f})$ for all $e\in \bar{E}(\Gamma)$ and $f\in \bar{E}(\Lambda)$. Then, using the decomposition in Lemma \ref{edge_decomposition}, we can say something about $\abs{E(\Theta^s)}$.

\begin{proposition}
\label{sub_same}
Let $e\in \bar{E}(\Gamma)$ and $f\in \bar{E}(\Lambda)$, then:
\begin{align*}
\abs{E(\Theta^s_{\sub}(e, f))} &\leq 100\cdot C_{\Gamma}\cdot C_{\Lambda}\cdot \abs{e}, \\
\abs{E(\Theta^s_{\super}(e, f))} &\leq 100\cdot C_{\Gamma}\cdot C_{\Lambda}\cdot \abs{f}.
\end{align*}
\end{proposition}

\begin{proof}
The two inequalities are symmetric so we just prove the first. For each component $g\subset \Theta^s_{\sub}(e, f)$, there corresponds an element in $\sub(e, f)$ in the following way: $\bar{o}(p_{\Lambda}\circ g) = (e, i)$ where $i\in \sub(e, f)$. By Lemma \ref{subwords}:
\[
\sub(e, f) = \bigcup_{i=1}^k\ap(p_i, q_i, r_i)
\]
where $k\leq \left\lfloor \frac{\abs{e}}{\abs{f}}\right\rfloor$ and $p_i\cdot q_i\leq \abs{f}$. Moreover, if $p_i\geq 2$, then $q = \per(f)$. But by Proposition \ref{same_period}, there can be at most $100\cdot C_{\Gamma}\cdot C_{\Lambda}$ elements in each arithmetic progression that correspond to components of $\Theta^s_{\sub}(e, f)$. Hence we get:
\[
\abs{E(\Theta^s_{\sub}(e, f))}\leq 100\cdot C_{\Gamma}\cdot C_{\Lambda}\cdot \left\lfloor \frac{\abs{f}}{\abs{e}}\right\rfloor\cdot \abs{e}\leq 100\cdot C_{\Gamma}\cdot C_{\Lambda}\cdot \abs{f}.
\]
\end{proof}

\begin{proposition}
\label{suff_same}
Let $e\in \bar{E}(\Gamma)$ and $f\in \bar{E}(\Lambda)$, then:
\[
\abs{E(\Theta^s_{\pref}(e, f))}, \abs{E(\Theta^s_{\suff}(e, f))} \leq 200\cdot C_{\Gamma~}\cdot C_{\Lambda}\cdot \min{\{\abs{e}, \abs{f}\}}.
\]
\end{proposition}

\begin{proof}
The two cases are symmetric so we just prove the inequality for $\Theta^s_{\suff}(e, f)$. For each component $g\subset \Theta^s_{\suff}(e, f)$, there corresponds an element $\suff(e, f)$ in the following way: $\bar{o}(p_{\Lambda}\circ g) = (e, i)$ where $i\in \suff(e, f)$. By Proposition \ref{suffix_overlap}:
\[
\suff(e, f) = \bigcup_{i = 1}^k\ap(p_i, q_i, r_i)
\]
where $k\leq \min{\{\lfloor\log_2(\abs{e})\rfloor, \lfloor\log_2(\abs{f})\rfloor\}}$, $q_i\leq \frac{q_{i-1}}{2}$ and $r_i \geq \abs{f} - q_{i-1}$. By Proposition \ref{same_period}, there can be at most $100\cdot C_{\Gamma}\cdot C_{\Lambda}$ many integers in each arithmetic progression that correspond to components $g\subset \Theta^s_{\suff}(e, f)$. Since $q_i\leq \frac{q_{i-1}}{2}$ and $r_i \geq \abs{f} - q_{i-1}$, we get:
\begin{align*}
\abs{E(\Theta^s_{\suff}(e, f))} &\leq 100\cdot C_{\Gamma}\cdot C_{\Lambda}\cdot \left(\min{\{\abs{e}, \abs{f}\}} + \frac{\min{\{\abs{e}, \abs{f}\}}}{2} + ... + \frac{\min{\{\abs{e}, \abs{f}\}}}{2^{\min{\{\lfloor\log_2(\abs{e})\rfloor, \lfloor\log_2(\abs{f})\rfloor\}}}}\right)\\
&\leq 200\cdot C_{\Gamma}\cdot C_{\Lambda}\cdot \min{\{\abs{e}, \abs{f}\}}.
\end{align*}
\end{proof}

\begin{proposition}
\label{boundary_edges}
\[
\sum_{e\in \bar{E}(\Gamma)}\sum_{f\in \bar{E}(\Lambda)}\abs{\Theta_{\partial}(e, f)}\leq 8\cdot C_{\Gamma}\cdot C_{\Lambda}\cdot \min{\{\abs{E(\Gamma)}, \abs{E(\Lambda)}\}}
\]
\end{proposition}

\begin{proof}
Each element $g\in \Theta_{\partial}(e, f)$ satisfies $\abs{g}\leq \min{\{\abs{e}, \abs{f}\}}$ by definition. By Lemma \ref{ve_bounds}, there are at most $8\cdot C_{\Gamma}\cdot C_{\Lambda}$ such elements ranging over all pairs $e\in \bar{E}(\Gamma)$, $f\in \bar{E}(\Lambda)$.
\end{proof}

Putting everything from this section together, we may prove one of our main theorems.

\begin{theorem}
\label{main_structure}
\[
\abs{E(\Theta^s)} \leq 508\cdot C_{\Gamma}\cdot C_{\Lambda}\cdot \left(\abs{\bar{E}(\Lambda)}\cdot \abs{E(\Gamma)} + \abs{\bar{E}(\Gamma)}\cdot \abs{E(\Lambda)}\right).
\]
\end{theorem}

\begin{proof}
By Propositions \ref{sub_same} and \ref{suff_same}:
\[
\abs{E(\Theta^s(e, f)) - E(\Theta_{\partial}(e, f))}\leq 500\cdot C_{\Gamma}\cdot C_{\Lambda}\cdot (\abs{e} + \abs{f}).
\]
Combining this with Proposition \ref{boundary_edges}, Lemma \ref{edge_decomposition} and the following:
\begin{align*}
\sum_{e\in \bar{E}(\Gamma)}\sum_{f\in \bar{E}(\Lambda)} (\abs{e} + \abs{f}) & = \sum_{e\in \bar{E}(\Gamma)} \left(\abs{\bar{E}(\Lambda)}\cdot \abs{e} + \abs{E(\Lambda)}\right)\\
													       & = \abs{\bar{E}(\Lambda)}\cdot \abs{E(\Gamma)} + \abs{\bar{E}(\Gamma)}\cdot \abs{E(\Lambda)},
\end{align*}
the result is proven.
\end{proof}

\subsection{Bounding the number of cycles}

Topological methods have been used for a long time to study fibre products of immersions of undirected graphs (see \cite{sta_83,neumann_90,tardos_92,dicks_01,khan_02,mineyev_12}).  By the Hanna Neumann Theorem, if $\gamma$ and $\lambda$ are immersions of undirected graphs, then $\core(\Gamma\times_{\Delta}\Lambda)$ has at most $\chi(\Gamma)\cdot \chi(\Lambda)$ many non cyclic components. However, handling the cyclic components requires a different approach. Let us first consider the simplest case when $\Gamma, \Lambda, \Delta\isom S^1$. Suppose $\deg(\gamma) = n\cdot p$ and $\deg(\lambda) = n\cdot q$ for integers $n, p, q\geq 1$ with $\gcd(p, q) = 1$. Then $\Theta = \Gamma\times_{\Delta}\Lambda$ has precisely $n$ components. So the number of components in general of the fibre product cannot depend on the topology of the graphs as $n$ was arbitrary. However, if $S_1, S_2\subset \Theta$ are any two distinct components, then 
\[
[\gamma\circ p_{\Gamma}|_{S_1}] = [\gamma\circ p_{\Gamma}|_{S_2}].
\]
Hence, all of the components in $\Theta$ represent the same conjugacy class of subgroup in $\pi_1(\Delta)$. We record a slight extension of this observation with the following Lemma.

\begin{lemma}
\label{rank_one_case}
If $\Gamma\isom \Lambda\isom S^1$, then the canonical map $\theta:\core(\Gamma\times_{\Delta}\Lambda)\to \Delta$ is the disjoint union of at most $\gcd(\deg(\gamma), \deg(\lambda))$ many copies of a single free homotopy class of map $S^1\to \Delta$.

\end{lemma}

This simpler case suggests that the better quantity to count is the number of distinct homotopy classes of maps (of connected graphs) that $\theta$ decomposes into. Indeed this line of approach is suggested in \cite{neumann_90} in the context of free groups. The aim of this section is to provide general bounds for this quantity. The main non-trivial step is to bound the number of cycles in $\Theta^s$. The key observation needed is the following: each element of $\bar{E}_c(\Theta^s)$ must support some lift in $\lift(p_{\Gamma}^s, e)$ for some $e\in \bar{E}(\Gamma)$ such that $p_{\Lambda}^s\circ\tilde{e}\in \olift(\lambda, \gamma\circ e)$.

\begin{proposition}
\label{overlaps}
Let $e\in \bar{E}(\Gamma)$, then
\[
\abs{\{\tilde{e}\in \lift(p_{\Gamma}^s, e) \mid p_{\Lambda}^s\circ\tilde{e}\in \olift(\lambda, \gamma\circ e)\}}\leq 500\cdot C_{\Gamma}\cdot C_{\Lambda}\cdot \abs{\bar{E}(\Lambda)}.
\]
\end{proposition}

\begin{proof}
If $\tilde{e}\in \lift(p^s_{\Gamma}, e)$ such that $p_{\Lambda}^s\circ\tilde{e}\in \olift(\lambda, \gamma\circ e)$, then 
\[
E(\tilde{e})\subset \bigcup_{f\in \bar{E}(\Lambda)} (E(\Theta^s_{\sub}(e, f)) \cup E(\Theta^s_{\pref}(e, f))\cup E(\Theta^s_{\suff}(e, f))).
\]
By Propositions \ref{sub_same} and \ref{suff_same}:
\[
\abs{E(\Theta^s_{\sub}(e, f)) \cup E(\Theta^s_{\pref}(e, f))\cup E(\Theta^s_{\suff}(e, f))}\leq 500\cdot C_{\Gamma~}\cdot C_{\Lambda}\cdot \abs{e},
\]
and so:
\begin{align*}
\bigcup_{\tilde{e}}E(\tilde{e}) &\leq 500\cdot \sum_{f\in \bar{E}(\Lambda)}500\cdot C_{\Gamma}\cdot C_{\Lambda}\cdot \abs{e}\\
					       & \leq 500\cdot C_{\Gamma}\cdot C_{\Lambda}\cdot \abs{\bar{E}(\Lambda)}\cdot \abs{e}.
\end{align*}
where $\tilde{e}$ ranges over the elements in the set $\{\tilde{e}\in \lift(p_{\Gamma}^s, e) \mid p_{\Lambda}^s\circ\tilde{e}\in \olift(\lambda, \gamma\circ e)\}$. But now since $E(\tilde{e}_1)\cap E(\tilde{e}_2) = \emptyset$ for every pair of distinct lifts $\tilde{e}_1, \tilde{e}_2\in \lift(p_{\Gamma}^s, e)$, the claimed bound follows.
\end{proof}

\begin{theorem}
\label{covering_overlaps}
There is a collection of paths $\{g_i:I\to \Theta^s\}_{i=1}^n$ such that 
\begin{enumerate}
\item $\cup_{i=1}^n g_i(I) = \Theta^s$, 
\item either $p_{\Gamma}\circ g_i$ factors through some $e\in \bar{E}(\Gamma)$ or $p_{\Lambda}\circ g_i$ factors through some $f\in \bar{E}(\Lambda)$,
\item $n\leq 1008\cdot C_{\Gamma}\cdot C_{\Lambda}\cdot \abs{\bar{E}(\Gamma)}\cdot \abs{\bar{E}(\Lambda)}$.
\end{enumerate}
\end{theorem}

\begin{proof}
Let $e\in E(\Theta^s)$ be any edge. Let $g_e, h_e:I\to \Theta$ be the largest paths containing $e$ and such that $p^s_{\Gamma}\circ g_e$ factors through some element in $\bar{E}(\Gamma)$ and $p^s_{\Lambda}\circ h_e$ factors through some element in $\bar{E}(\Lambda)$. We can see that if $p^s_{\Gamma}\circ g_e$ is not equal to an element in $\bar{E}(\Gamma)$, then $g$ has an end point on the boundary of $\Theta$. Similarly for $h_e$. So let us decompose $E(\Theta^s) = E_1\sqcup E_2$ where $E_1$ contains all the edges $e$ such that either $g_e$ or $h_e$ has an endpoint on the boundary. By Lemma \ref{ve_bounds}, we can cover all edges in $E_1$ by at most $8\cdot C_{\Gamma}\cdot C_{\Lambda}$ of the required paths. 

Now for every $e\in E_2$, we have that $p^s_{\Gamma}\circ g_e\in \bar{E}(\Gamma)$ and $p^s_{\Lambda}\circ h_e\in\bar{E}(\Lambda)$. So if $\abs{g_e}\geq \abs{h_e}$, then the projection of $g_e$ to $\Lambda$ must traverse a vertex in $\bar{V}(\Lambda)$. In other words:
\[
g_e\in \{\tilde{e}\in \lift(p_{\Gamma}^s, p^s_{\Gamma}\circ g_e) \mid p_{\Lambda}^s\circ\tilde{e}\in \olift(\lambda, \gamma\circ p^s_{\Gamma}\circ g_e)\}.
\]
If $\abs{g_e}\leq \abs{h_e}$ we have:
\[
h_e\in \{\tilde{e}\in \lift(p_{\Lambda}^s, p^s_{\Lambda}\circ h_e) \mid p_{\Gamma}^s\circ\tilde{e}\in \olift(\gamma, \lambda\circ p^s_{\Lambda}\circ h_e)\}.
\]
We may bound the sets on the right by Proposition \ref{overlaps}. Now we take the set $\{g_i:I\to \Theta\}_{i=1}^n$ to consist of all distinct paths in the following set:
\[
\bigcup_{e\in \Theta^s}\left(\left(\bigcup_{\abs{g_e}\geq\abs{h_e}}g_e\right)\cup \left(\bigcup_{\abs{g_e}<\abs{h_e}}h_e\right)\right).
\]
Summing over elements of $\bar{E}(\Gamma)$ and $\bar{E}(\Lambda)$ we get:
\begin{align*}
n \leq & \sum_{e\in \bar{E}(\Gamma)} \abs{\{\tilde{e}\in \lift(p_{\Gamma}^s, e) \mid p_{\Lambda}^s\circ\tilde{e}\in \olift(\lambda, \gamma\circ e)\}}\\
&+ \sum_{f\in \bar{E}(\Lambda)} \abs{\{\tilde{f}\in \lift(p_{\Lambda}^s, f) \mid p_{\Gamma}^s\circ\tilde{f}\in \olift(\gamma, \lambda\circ f)\}}+ 8\cdot C_{\Gamma}\cdot C_{\Lambda}\\
& \leq 1008\cdot C_{\Gamma}\cdot C_{\Lambda}\cdot \abs{\bar{E}(\Gamma)}\cdot \abs{\bar{E}(\Lambda)}.
\end{align*}
\end{proof}

A simplified version of the proof of Theorem \ref{covering_overlaps} yields the following bound on the number of cycles in $\Theta^s$.

\begin{theorem}
\label{cycles}
\[
\abs{\bar{E}_c(\Theta^s)}\leq 500\cdot \left(\abs{\bar{E}(\Gamma)}\cdot \abs{\bar{E}(\Lambda)}\right)^2.
\]
\end{theorem}

\begin{proof}
Without loss, we may assume that $V^m(\Gamma), V^m(\Lambda) = \emptyset$, thus we have $C_{\Gamma} = \abs{\bar{E}(\Gamma)}$ and $C_{\Lambda} = \abs{\bar{E}(\Lambda)}$. Each element $g\in \bar{E}_c(\Theta^s)$ must support some element in $\{\tilde{e}\in \lift(p_{\Gamma}^s, e) \mid p_{\Lambda}^s\circ\tilde{e}\in \olift(\lambda, \gamma\circ e)\}$ for some $e\in \bar{E}(\Gamma)$. Hence, just as in the proof of Theorem \ref{covering_overlaps}, we sum over elements of $\bar{E}(\Gamma)$ to complete the proof.
\end{proof}

\begin{theorem}
\label{full_homotopy_components}
The natural map $\theta:\Theta\to \Delta$ contains at most
\[
40538\cdot (\beta_1(\Gamma)\cdot \beta_1(\Lambda))^2
\]
many free homotopy classes of maps of connected graphs.
\end{theorem}

\begin{proof}
Contractible components all lie in the same free homotopy class, so we add one to our count and assume that $V^m(\Gamma), V^m(\Lambda) = \emptyset$. By Lemma \ref{betti_bound}, we have $\abs{\bar{E}(\Gamma)}\leq 3\cdot\beta_1(\Gamma)$ and $\abs{\bar{E}(\Lambda)}\leq 3\cdot \beta_1(\Lambda)$. By Lemma \ref{ve_bounds}, we have $\abs{\bar{E}(\Theta) - \bar{E}_c(\Theta)}\leq 36\cdot \beta_1(\Gamma)\cdot \beta_1(\Lambda)$. So now we simply need to bound the free homotopy classes of maps corresponding to the cycles $\bar{E}_c(\Theta)$. By Lemma \ref{rank_one_case}, and Theorem \ref{long_cycles_graph}, the set $\{\gamma\circ p_{\Gamma}\circ f \mid f\in\bar{E}_c(\Theta)\}$ contains at most $\beta_1(\Gamma)\cdot \beta_1(\Lambda)$ many free homotopy classes of cycles that factor through $s(\mathbb{S})$. Now applying Theorem \ref{cycles} finishes the proof.
\end{proof}

We conclude this section with an example demonstrating a lower bound for a sharp upper bound of the free homotopy classes of cycles in $\theta:\Theta\to \Delta$.

\begin{example}
\label{lower_bound}
Let $k, m\geq 1$ be a pair of integers. Let $\Delta$ be the connected graph with a single vertex and $k + m$ edges. We may label these loops by $x_1, ..., x_k, y_1, ..., y_m$. Consider the paths: 
\begin{align*}
w &= x_1* ...* x_k* y_1*...* y_m\\
w_i &= w[i:]* w[:i].
\end{align*}
Let $p_1, ..., p_{k+m}$ be a collection of distinct primes. Consider the connected graph $\Gamma$ with $\bar{V}(\Gamma) = \{v_{\Gamma}\}$ and $\bar{E}(\Gamma) = \{w_1^{p_1}, ..., w_k^{p_k}\}$. Consider also the connected graph $\Lambda$ with $\bar{V}(\Lambda) = \{v_{\Lambda}\}$ and $\bar{E}(\Lambda) = \{w_{k+1}^{p_{k+1}}, ..., w_{k+m}^{p_{k+m}}\}$. Then $\rcore(\Gamma\times_{\Sigma}\Lambda)$ has precisely 
\[
k\cdot m = \beta_1(\Gamma)\cdot \beta_1(\Lambda)
\]
components. By choosing an appropriate basepoint for each component of $\rcore(\Gamma\times_{\Delta}\Lambda)$, the labels of each cycle are $w^{p_i\cdot p_j}$ as $i$ ranges from $1$ to $k$ and $j$ ranges from $k+1$ to $k+m$. Hence, the induced maps to $\Delta$ restricted to each component are in distinct free homotopy classes. Furthermore, each component represents conjugates $\pi_1(\Gamma)\cap \pi_1(\Lambda)^g$, none of which are conjugate into each other. On the other hand, they are all conjugate into the cyclic group $\langle w\rangle$.
\end{example}

\section{Free groups}
\label{free_groups_section}

Now that we have proved all of our results in the context of graphs, we can discuss applications of these results to free groups. We begin by making the connection with free groups explicit. An \emph{involutive graph} $\Gamma$ is a graph, together with an involution $^{-1}:E(\Gamma)\to E(\Gamma)$ such that $o(e) = t(e^{-1})$ and $e^{-1} \neq e$ for all $e\in E(\Gamma)$. An involutive graph map $\gamma:\Gamma\to \Delta$ must also satisfy
\[
\gamma(e^{-1}) = \gamma(e)^{-1}
\]
for all $e\in E(\Gamma)$.

A path $f:I\to \Gamma$ is \emph{reduced} if there is no subsequent pair of edges $e_1, e_2\in E(I)$ such that $f(e_1) = f(e_2)^{-1}$. The core of an involutive graph is instead defined as the union of all the reduced cycles. The relative core is defined similarly. The set of (endpoint preserving homotopy classes of) reduced paths forms a groupoid, the fundamental groupoid $\pi_1(\Gamma)$. If we choose a basepoint $x\in V(\Gamma)$, then the set of (basepoint preserving homotopy classes of) reduced paths actually form a group, the fundamental group $\pi_1(\Gamma, x)$. By Stallings construction \cite{sta_83}, there is a bijection:
\[
\{\gamma:(\Gamma, x)\to (\Delta, y) \mid \gamma \text{ based immersion of relatively core involutive graphs}\} \to \{ A < \pi_1(\Delta, y)\}
\]
given by the $\pi_1$ functor. And so, as elucidated in detail in \cite{kapovich_02}, results for involutive graphs correspond directly to results for subgroups of free groups. 

Since the graphs we have worked with were not involutive graphs, we will instead consider a slight modification of the classic construction. So let $\Gamma$ be a graph equipped with an involution on its edges $^{-1}:E(\Gamma)\to E(\Gamma)$  and its vertices $^{-1}:V(\Gamma)\to V(\Gamma)$ such that
\begin{enumerate}
\item\label{itm:inverse_edge_vertex} for every $e\in E(\Gamma)$, we have $o(e) = t(e^{-1})^{-1}$,
\item\label{itm:inverse_edge} for all $e\in E(\Gamma)$, $(e^{-1})^{-1} = e$ and $e^{-1}\neq e$,
\item\label{itm:inverse_vertex} for every $v\in V(\Gamma)$, $(v^{-1})^{-1} = v$.
\end{enumerate}
We will call this a \emph{fully involutive graph}. Note that by \ref{itm:inverse_edge_vertex}, the involution $^{-1}$ can be extended to the set of paths: if $f:I\to \Gamma$ is a path, then $f^{-1}$ is also a path and $o(f) = t(f^{-1})$.

Let $\Gamma$ be a fully involutive graph. We may also define reduced paths in $\Gamma$ in the exact same way as we did for involutive graphs. However, the fundamental groupoid may not always be defined. If $x\in V(\Gamma)$ such that $x^{-1} = x$, then by \ref{itm:inverse_edge_vertex}, the fundamental group $\pi_1(\Gamma, x)$ is well defined. 

There is a standard formula relating the rank of the fundamental group of an involutive graph to its combinatorial data. We record this formula also for fully involutive graphs.

\begin{lemma}
\label{easy_rank}
Let $\Gamma$ be a connected fully involutive graph. Then:
\[
\rk(\Gamma, x) = \frac{\abs{\bar{E}(\Gamma) - \bar{E}_c(\Gamma)}}{2} - \abs{\bar{V}(\Gamma)} + 1.
\]

\end{lemma}

The quotient map $\Gamma\to \Gamma'$ which identifies all $v\in V(\Gamma)$ with $v^{-1}$, induces a functor from the category of fully involutive graphs to the category of involutive graphs. Denote this functor by $ST$.

\begin{figure}
\centering
\vspace*{8pt}
\includegraphics[scale=0.8]{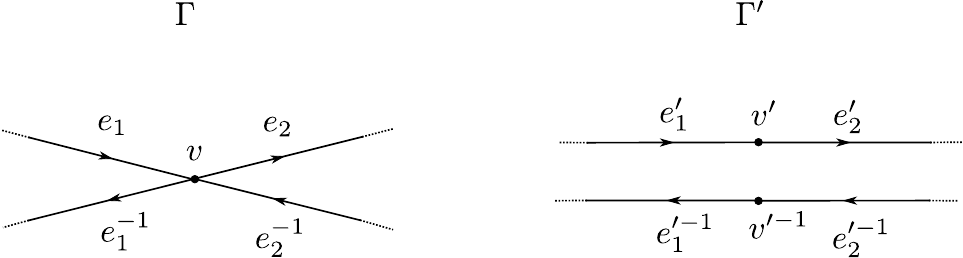}
\vspace{8pt}
\caption{A vertex splitting}
\label{vertex_splitting}
\end{figure}

We say a fully involutive graph $\Gamma'$ is obtained from a fully involutive graph $\Gamma$ by \emph{splitting a vertex} if: there is some vertex $v\in V(\Gamma) - V^m(\Gamma)$ such that $\bstar(v) = \{e_1, e_2^{-1}\}$ and $\fstar(v) = \{e_1^{-1}, e_2\}$ and $\Gamma'$ is obtained from $\Gamma$ by removing the vertex $v$, removing the edges $e_1$ and $e_2$, adding two new vertices $v'$, $v'^{-1}$ and four new edges $e_1', e_1'^{-1}, e_2'$ and $e_2'^{-1}$ such that:
\begin{align*}
o(e_1') &= o(e_1), & t(e_1') &= v',\\
o(e_1'^{-1}) &= v'^{-1}, & t(e_1'^{-1}) &= t(e_1^{-1}),\\
o(e_2') &= v', & t(e_2') &= t(e_2),\\
o(e_2'^{-1}) &= o(e_2^{-1}), & t(e_2'^{-1}) &= v'^{-1}.
\end{align*}
This operation is perhaps better understood through a diagram, see Figure \ref{vertex_splitting}. There is an induced quotient map $q:\Gamma'\to \Gamma$ where $v'$ and $v'^{-1}$ map to $v$, $e_1', e_1'^{-1}, e_2', e_2'^{-1}$ map to $e_1, e_1^{-1}, e_2, e_2^{-1}$ respectively and everything else maps where expected. The induced quotient map will be called a \emph{vertex splitting}. A graph with involution $\Gamma$ is \emph{unsplittable} if no vertex splitting can be performed on $\Gamma$.

\begin{lemma}
\label{unique_splitting}
Let $\Gamma$ be a fully involutive graph. There is a sequence of vertex splittings:
\[
\Gamma_n\to ...\to \Gamma_1\to \Gamma.
\]
such that $\Gamma_n$ is unsplittable. Moreover, $\Gamma_n$ is uniquely defined.
\end{lemma}

\begin{proof}
Since each new vertex created in a vertex splitting cannot be split again, there is a sequence of vertex splittings $\Gamma_n\to ...\to \Gamma_1\to \Gamma$ such that $\Gamma_n$ is unsplittable and $n\leq \abs{V(\Gamma)}$. Uniqueness follows from the fact that any pair of vertex splittings commute.
\end{proof}

\begin{figure}
\centering
\vspace*{8pt}
\includegraphics[scale=0.8]{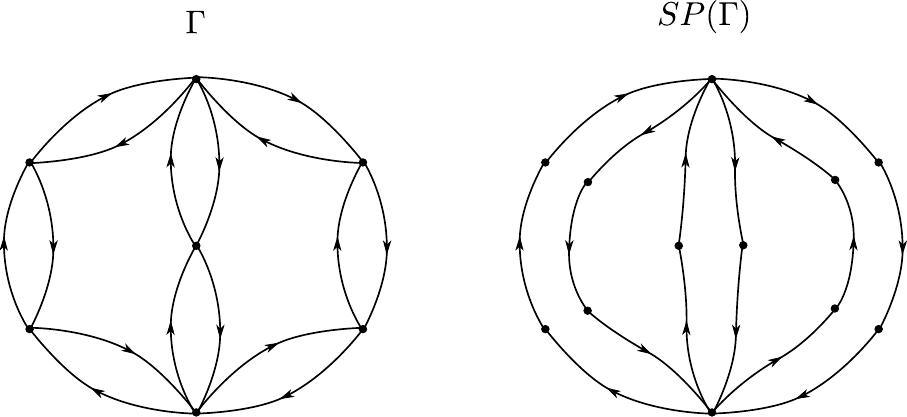}
\vspace{8pt}
\caption{Applying the functor $SP$}
\label{splitting_stallings}
\end{figure}

By Lemma \ref{unique_splitting}, we will write $SP(\Gamma)$ for the unique unsplittable fully involutive graph obtained from a fully involutive graph $\Gamma$. See Figure \ref{splitting_stallings} for an example of this. This is also a functor and a section for the functor $ST$. We record this in the following Lemma.

\begin{lemma}
\label{functoriality}
$SP$ and $ST$ are functors, with $SP$ a section for $ST$.

\end{lemma}

We are now ready to prove our main Theorems for subgroups of free groups from the introduction.

\begin{theorem}
\label{finite_index}
Let $F$ be a free group and $A<F$ a finitely generated subgroup. Then there are at most $\rk(A)$ many conjugacy classes of maximal cyclic subgroups $M< F$ such that 
\[
12\cdot\rk(A)\leq[M:A\cap M]<\infty.
\]
\end{theorem}

\begin{proof}
Let $\gamma:(\Gamma, y)\to (\Delta, x)$ be the relatively core involutive graph immersion realising the inclusion $A<F$. Let $[G_1], [G_2], ... $ be conjugacy classes of maximal cyclic subgroups such that $12\cdot \rk(A)\leq[G_i:A\cap G_i]<\infty$. Now let $g_i$ be a generator for $G_i$. By construction no pair $g_i$, $g_j$ are conjugate to each other or to an inverse. Furthermore, $o_A(g)\geq12\cdot \rk(A)$ where $o_A(g)$ is the relative order of $g$ defined in the introduction. Let $h_i$ be a shortest representative in the conjugacy class $[g_i]$ in $F$. Then there are reduced cycles $f_1, ..., f_n:S^1\to SP(\Gamma)$ such that (after choosing an appropriate basepoint for $S^1$) $[f_i] = h_i$.

Let $w_i:I\to SP(\Delta)$ be a primitive loop such that $f_i$ is a $w_i$-cycle. Assume that $\abs{w_i}\geq \abs{w_{i+1}}$. By Lemma \ref{long_w_edge}, for each $f_i$ there is an element $e_i\in \bar{E}(SP(\Gamma))$ such that $\abs{e_i}\geq 2\cdot \abs{w_i}$ and $f_i$ traverses $e_i$. If $f_{j}$ traverses $e_i$ for any $i<j$, up to replacing $w_i$ with a cyclic conjugate of itself, $w_i^2$ has period $\abs{w_j}\leq \abs{w_i}$. Since $w_i$ and $w_j$ are primitive, it follows from Corollary \ref{square_period} that $w_i$ is a cyclic conjugate of $w_j$. But since $\gamma$ was an immersion, it then also follows that $f_i = f_j$. Similarly if $f_j$ traverses $e_i^{-1}$, but the contradiction is that $f_i = f_j^{-1}$. Hence, $f_i$ is supported in $SP(\Gamma) - \{e_1, e_1^{-1}, ..., e_{i-1}, e_{i-1}^{-1}\}$. As $SP(\Gamma) - \{e_1, e_1^{-1}, ..., e_{\rk(A)}, e_{\rk(A)}^{-1}\}$ cannot support any reduced non-trivial cycles, the result follows.
\end{proof}

\begin{corollary}
\label{relative_spectrum}
Let $F$ be a free group and $A<F$ a finitely generated subgroup. Then there exists a set $N\subset \N$ of cardinality bounded above by $13\cdot \rk(A)$, such that:
\[
\mathcal{O}_{A}(F) = \bigcup_{n\in N}\bigcup_{d\mid n} d.
\]

\end{corollary}

\begin{theorem}
\label{intersection_rank}
There exists a universal constant $C>0$, with the following property. For any free group $F$ and any finitely generated subgroups $A, B<F$, we have:
\[
\sum_{[A\cap B^g]} \rk(A\cap B^g) \leq C\cdot (\rk(A)\cdot \rk(B))^2.
\]
\end{theorem}

\begin{proof}
Let $\Delta$ be an involutive graph with a single vertex and such that $\pi_1(\Delta, x)\isom F$. Let $\gamma:(\Gamma, y)\to (\Delta, x)$ and $\lambda:(\Lambda, z)\to (\Delta, x)$ be the relatively core involutive graph immersions realising the inclusions of the subgroups $A$ and $B<\pi_1(\Delta, x)$ respectively. By Theorem \ref{full_homotopy_components}, the map:
\[
\theta: \rcore(SP(\Gamma)\times_{SP(\Delta)}SP(\Lambda))\to SP(\Delta)
\]
consists of at most $C\cdot (\beta_1(SP(\Gamma))\cdot \beta_1(SP(\Lambda)))^2$ many free homotopy classes of connected graph maps where $C>0$ is a universal constant. Now by Lemma \ref{functoriality}, we have that the map:
\[
\rcore(\Gamma\times_{\Delta}\Lambda)\to \Delta
\]
is the same as the map $ST(\theta)$. Hence it must also consist of at most $C\cdot (\beta_1(SP(\Gamma))\cdot \beta_1(SP(\Lambda)))^2$ many free homotopy classes of connected graph maps. Thus by Theorem 5.5 in \cite{sta_83} and Lemma \ref{easy_rank}, the Theorem follows.
\end{proof}

Effectively what we have done in this section is convert our statements of our main fibre product results to hold also for immersions of undirected graphs. Therefore, we record this in the following Theorem.

\begin{theorem}
Theorem \ref{long_cycles_graph_equiv} also holds for graph immersions of undirected graphs. Theorems \ref{main_structure} and \ref{full_homotopy_components} also hold for graph immersions of undirected graphs, possibly with larger constants.
\end{theorem}

\section{algorithmic results}
\label{algorithm_section}

In this section we present our algorithmic results. In what follows, all algorithms will have as input two forwards graph immersions $\gamma:\Gamma\to \Delta$ and $\lambda:\Lambda\to \Delta$. We will also use the same notation given in the beginning of Section \ref{main_structure_section}. The parameters for our complexity results will be the following two quantities:
\begin{align*}
m &= (\beta_1(\Gamma) + \abs{\partial\Gamma}) + (\beta_1(\Lambda) + \abs{\partial\Lambda}),\\
n &= \abs{E(\Gamma)} + \abs{E(\Lambda)}.
\end{align*}

\subsection{Accepting paths for $w$-cycles}

Let $w:I\to \Delta$ be a primitive loop. Recall that we denote by $\hat{w}:S^1\to \Delta$ the graph map obtained from $w$ by identifying the vertices in $\partial I$ and a $w$-cycle is a graph map $f:S^1\to \Delta$ that factors through $\hat{w}:S^1\to \Delta$. An important part of our DFA-NEI algorithm is the following variation of Theorem 2 in \cite{oliveira_18}.

\begin{theorem}
\label{unary_generalisation}
Let $w:I\to \Delta$ be a primitive loop and $\gamma:\Gamma\to \Delta$ and $\lambda:\Lambda\to \Delta$ be $w$-cycles. If all the initial states in both $\Gamma$ and $\Lambda$ are contained in at most $k$ segments of length $\abs{w}$, then we may decide in $O(k^2\cdot n)$ time if $\core(\Gamma\times_{\Delta}\Lambda)$ contains an accepting path.
\end{theorem}

\begin{proof}
We may find shortest loops $w_{\Gamma}:I\to \Gamma$ and $w_{\Lambda}:I\to \Lambda$ such that $w_{\Gamma}[:\abs{w}] = w_{\Lambda}[:\abs{w}]$ in $O(n)$ time by Proposition \ref{cyclic_conjugation_algorithm}. We may associate to each vertex $v\in V(\Gamma)$ a pair of integers $(i, j)$ in the following way: if $0\leq l<\abs{w_{\Gamma}}$ is the unique integer such that $w_{\Gamma}(l) = v$, then let $i = \left\lfloor\frac{l}{\abs{w}}\right\rfloor$ and $j = l - i$. Denote by $\eta_{\Gamma}:V(\Gamma)\to \N\times\N$ the corresponding map. Similarly we have a map $\eta_{\Lambda}:V(\Lambda)\to \N\times\N$ for $\Lambda$. By assumption, we have $\abs{p_1\circ\eta_{\Gamma}(V^i(\Gamma))}\leq 2\cdot k$ and $\abs{p_1\circ\eta_{\Lambda}(V^i(\Lambda))}\leq 2\cdot k$ where $p_1:\N\times\N\to \N$ is projection to the first factor. Both the sets $p_1\circ\eta_{\Gamma}(V^i(\Gamma))$ and $p_1\circ\eta_{\Lambda}(V^i(\Lambda))$ may also be computed in $O(n)$ time.

For each pair $a\in p_1\circ\eta_{\Gamma}(V^i(\Gamma))$, $b\in p_1\circ\eta_{\Lambda}(V^i(\Lambda))$, we check if there exists some $0\leq c< \abs{w}$ such that $(a,c)\in \eta_{\Gamma}(V^i(\Gamma))$ and $(b, c)\in \eta_{\Lambda}(V^i(\Lambda))$. This may be done in $O(n)$ time for each choice of $a, b$. If there is no such pair, then $V^i(\core(\Gamma\times_{\Delta}\Lambda)) = \emptyset$ and $\core(\Gamma\times_{\Delta}\Lambda)$ cannot contain an accepting path.

Now let $a$, $b$ be a pair of such integers. With some minor modifications, we now essentially follow the proof of Theorem 2 in \cite{oliveira_18}. Let $\gamma\circ w_{\Gamma} = w^m$ and $\lambda\circ w_{\Lambda} = w^n$ and $q = \gcd(m, n)$. Then $\core(\Gamma\times_{\Delta}\Lambda)$ contains an accepting path starting from some vertex in $V(w_{\Gamma}[a\cdot \abs{w}:(a+1)\cdot \abs{w} -1])\times V(w_{\Lambda}[b\cdot \abs{w}:(b+1)\cdot \abs{w} -1])$ if and only if there is a vertex $(x, y)\in V^f(\Gamma)\times V^f(\Lambda)$ and integers $s, t, r, c$ such that:
\begin{align*}
\eta_{\Gamma}(x) &= (s\cdot q + r - a, c),\\
\eta_{\Lambda}(y) &= (t\cdot q + r - b, c).
\end{align*}
Now for each $x\in V^f(\Gamma)$, we can compute the unique minimal pair $(r', c')$ such that $\eta_{\Gamma}(s'\cdot q + r' - a, c')$ for some $s'$. We may put all these integers in the set $S_{\Gamma}$. Similarly we compute the set $S_{\Lambda}$. We may do this in $O(n)$ time. Finally, we decide if $S_{\Gamma}\cap S_{\Lambda} = \emptyset$ or not. This may also be done in $O(n)$ time.

Thus, by iterating over all pairs $(a, b) \in p_1\circ\eta_{\Gamma}(V^i(\Gamma))\times p_1\circ\eta_{\Lambda}(V^i(\Lambda))$, we may decide if $\core(\Gamma\times_{\Delta}\Lambda)$ has any accepting paths in $O(k^2\cdot n)$ time.
\end{proof}

\subsection{Connectivity in $s(\mathbb{S})$}

As a product of the proof of Theorem \ref{long_cycles_graph}, we also obtain that the graph $s_{\Gamma}:\mathbb{S}_{\Gamma}\to \Gamma$ is computable. The following result states that we may also compute them efficiently.

\begin{theorem}
\label{long_cycles_algorithm}
There is an algorithm that computes the graph $s_{\Gamma}:\mathbb{S}_{\Gamma}\to \Gamma$ from Theorem \ref{long_cycles_graph} in $O(m^2\cdot n)$ time.
\end{theorem}

\begin{proof}
Since all cycles are supported in $\core(\Gamma)$, we may assume that $\Gamma$ is core itself. Let $f:S^1\to \Gamma$ be a cycle in $\mathbb{S}_{\Gamma}$. Then there is some primitive loop $w:I\to \Delta$ and some map $g:S^1\to S^1$ such that $\gamma\circ f = \hat{w}\circ g$ and $\deg(g)\geq12\cdot \beta_1(\Gamma)$. By Lemma \ref{long_w_edge}, there is some element $e\in \bar{E}(\Gamma)$ such that $e$ is a $w$-path and $\abs{e}\geq 2\cdot \abs{w}$. By Theorem \ref{periodicity}, it follows that $\per(e) = \abs{w}$. So in order to find the cycles in $\mathbb{S}_{\Gamma}$ we may do the following. By Theorem \ref{linear_pattern} we may compute $\per(e)$ for each $e\in \bar{E}(\Gamma)$ in $O(\abs{E(\Gamma)})$ time. Then let $w_e:I\to e$ denote a subpath of $e$ of length exactly $\per(e)$. By Lemma \ref{gamma_w_factorisation}, we may decide in $O\left(\abs{E(\Gamma)} + 2\cdot \abs{\bar{E}(\Gamma)}\cdot \abs{w_e}\right)$ time if this extends to a $w_e$-cycle. Since $\sum_{e\in \bar{E}(\Gamma)}\abs{w_e}\leq \abs{E(\Gamma)}$, we may do this for all $e\in \bar{E}(\Gamma)$ in $O\left(\abs{\bar{E}(\Gamma)}\cdot \abs{E(\Gamma)}\right)$ time. Moreover, we may also decide at the same time if the projections of the computed cycles to $\Delta$ have degree greater than or equal to $12\cdot \beta_1(\Gamma)$. Thus, we may decide if they are cycles in $\mathbb{S}_{\Gamma}$. Finally, by Proposition \ref{cyclic_conjugation_algorithm}, we may decide in $O\left(\abs{\bar{E}(\Gamma)}^2\cdot \abs{E(\Gamma)}\right)$ time if any of the computed cycles are the same. Hence overall, the algorithm will require $O\left(\abs{\bar{E}(\Gamma)}^2\cdot \abs{E(\Gamma)}\right)$ time to compute $s_{\Gamma}:\mathbb{S}_{\Gamma}\to \Gamma$. Combined with Lemma \ref{betti_bound} we obtain the result.
\end{proof}

The subgraph $s(\mathbb{S})\subset \Theta$ in general has more than a linear number of states and transitions, providing an obstacle to a linear solution to the DFA non-empty intersection problem. Therefore, in order to get around this, we may like to learn something about the graph $s(\mathbb{S})$ without actually computing it. We show that we may answer any question about connectivity of $s(\mathbb{S})$ efficiently.

\begin{proposition}
\label{long_graph_directions}
There is an integer $q> 0$ such that there is a $O(m^q\cdot n)$ time algorithm that, given as input elements $(e, f)\in E(\Gamma)\times E(\Lambda)$ and $(x, y)\in V(\Gamma)\times V(\Lambda)$, decides if there exists a path $g:I\to \Gamma\times_{\Delta}\Lambda$ with the following properties:
\begin{enumerate}
\item $g$ lifts to $s:\core(\mathbb{S}_{\Gamma}\times_{\Delta}\mathbb{S}_{\Lambda})\to \core(\Gamma\times_{\Delta}\Lambda)$, 
\item $g[0] = (e, f)$,
\item $t(g) = (x, y)$.
\end{enumerate}
\end{proposition}

\begin{proof}
Without loss we may assume that $\Gamma$ and $\Lambda$ are core graphs. By Theorem \ref{long_cycles_algorithm}, we may compute the graphs $s_{\Gamma}:\mathbb{S}_{\Gamma}\to \Gamma$ and $s_{\Lambda}:\mathbb{S}_{\Lambda}\to \Lambda$ in $O(m^2\cdot n)$ time. We set 
\begin{align*}
V^i(\mathbb{S}_{\Gamma}) &= \cup_{e'\in E(s_{\Gamma}^{-1}(e))}o(e'),\\
V^i(\mathbb{S}_{\Lambda}) &= \cup_{f'\in E(s_{\Lambda}^{-1}(f))}o(f'),\\
V^f(\mathbb{S}_{\Gamma}) &= s_{\Gamma}^{-1}(x),\\
V^f(\mathbb{S}_{\Lambda}) &= s_{\Lambda}^{-1}(y).
\end{align*}
Now by construction, the path $g$ exists if and only if there exists some admissible path in the graph $\mathbb{S} = \core(\mathbb{S}_{\Gamma}\times_{\Delta}\mathbb{S}_{\Lambda})$ with $V^i(\mathbb{S}) = V^i(\mathbb{S}_{\Gamma})\times V^i(\mathbb{S}_{\Lambda})$ and $V^f(\mathbb{S}) = V^f(\mathbb{S}_{\Gamma})\times V^f(\mathbb{S}_{\Lambda})$. We may decide which pairs of components of $\mathbb{S}_{\Gamma}$ and $\mathbb{S}_{\Lambda}$ factor through the same cycle in $O(m^2\cdot n)$ time using the algorithms from Theorem \ref{linear_pattern} and Proposition \ref{cyclic_conjugation_algorithm}. By Theorem \ref{extension_language}, if $S_{\Gamma}\subset \mathbb{S}_{\Gamma}$ and $S_{\Lambda}\subset \mathbb{S}_{\Lambda}$ are both $w$-cycles for some primitive loop $w:I\to \Delta$, there exist at most $10\cdot \left(\abs{\bar{E}(\Gamma)} + \abs{\bar{E}(\Lambda)}\right)\in O(m)$ many segments of length $\abs{w}$ that cover $V^i(S_{\Gamma})$ and $V^i(S_{\Lambda})$. Now the problem may be solved with the required complexity bounds by Theorem \ref{unary_generalisation}.
\end{proof}

\begin{proposition}
\label{adjacency_partition}
There is an integer $q> 0$ such that there is a $O(m^q\cdot n)$ time algorithm that computes subsets $V_{\mathbb{S}}\subset \bar{V}(\Gamma)\times \bar{V}(\Lambda)$ and $E_{\mathbb{S}}\subset E(\Gamma)\times E(\Lambda)$ along with partitions $V_{\mathbb{S}} = V_1\sqcup ...\sqcup V_c$ and $E_{\mathbb{S}} = E_1\sqcup...\sqcup E_c$ so that for every component $S\subset s(\mathbb{S})$ there is an $i$ such that:
\begin{enumerate}
\item $V_i = S\cap (\bar{V}(\Gamma)\times \bar{V}(\Lambda))$,
\item $E_i = S\cap (\bigcup_{v\in \bar{V}(\Gamma)\times\bar{V}(\Lambda)}\fstar(v, \Gamma\times_{\Delta}\Lambda))$.
\end{enumerate}
\end{proposition}

\begin{proof}
We may simply iterate the algorithm from Proposition \ref{long_graph_directions} for every pair of elements $u, v\in \bar{V}(\Gamma)\times \bar{V}(\Lambda)$ and every edge in $\fstar(u, \Gamma\times_{\Delta}\Lambda)$. 
\end{proof}

\subsection{Regular languages}

Let $\Sigma$ be a finite set. A \emph{finite state automaton over $\Sigma$} is a graph $\Gamma$ together with an initial vertex $V^i(\Gamma) = \{x\}$, a collection of final vertices $V^f(\Gamma)$ and a labelling of its edges by elements of $\Sigma$. If $\Delta$ is a graph with a single vertex and an edge for each element in $\Sigma$, then labellings of edges of $\Gamma$ by elements of $\Sigma$ are in a one-to-one correspondence with graph maps $\gamma:\Gamma\to \Delta$. Thus, for our purposes, a finite state automaton will be a quadruple $(\Gamma, \gamma, x, V^f)$, where $\gamma:\Gamma\to \Delta$ is a graph map with $\Delta$ a connected graph with a single vertex. A \emph{deterministic finite state automaton} (or \emph{DFA}) is an automaton $(\Gamma, \gamma, x, V^f)$ such that $\gamma$ is a forwards immersion. The \emph{language} of the automaton $(\Gamma, \gamma, x, V^f)$ is the following collection of paths:
\[
L(\Gamma) = \{\gamma\circ g \mid g:I\to \Gamma \text{ an accepting path}\}.
\]
It is a classical result of Rabin-Scott \cite{rabin_59} that, given two automata $(\Gamma, \gamma:\Gamma\to \Delta, x, V^f(\Gamma))$, $(\Lambda, \lambda:\Lambda\to \Delta, y, V^f(\Lambda))$, we have:
\[
L(\Gamma)\cap L(\Lambda) = L(\Gamma\times_{\Delta}\Lambda).
\]
If $\gamma$ and $\lambda$ are finite sheeted covers, then $\Gamma\times_{\Delta}\Lambda$ is a core graph and $V(\Gamma\times_{\Delta}\Lambda) = V(\Gamma)\cdot V(\Lambda)$. So the optimal complexity bound for computing the automaton representing $L(\Gamma)\cap L(\Lambda)$ is $O(n^2)$. 

On the other hand, it is an open problem if deciding if $L(\Gamma)\cap L(\Lambda)\neq \emptyset$ can be done in subquadratic time. See \cite{kara_00,kara_02,wehar_14,oliveira_18,oliveira_20} for results on this problem and its variations. Currently, the only known positive result for two DFAs is a linear time solution for automata over a unary alphabet, or rather, automata $\gamma:\Gamma\to \Delta$ where $\Delta$ has only a single edge. This constraint on the edges of $\Delta$ immediately puts constraints on the topology of $\Gamma$. Indeed, there are at most four distinct forms that $\Gamma$ can take: a point, a segment, a circle or a circle with a segment attached. Thus, the $m$ parameter is forced to be bounded above by $4$. Generalising the unary case, we show that for any fixed $m$, the DFA-NEI problem can be solved in linear time.

\begin{theorem}
\label{intersection_emptiness}
There is an integer $q>0$ such that there is a $O(m^q\cdot n)$ time algorithm that, given as input two deterministic finite state automata $\gamma:\Gamma\to \Delta$ and $\lambda:\Lambda\to \Delta$, decides if $L(\Gamma)\cap L(\Lambda)\neq \emptyset$.
\end{theorem}

\begin{proof}
Let us first fix some notation. As $\Gamma$ and $\Lambda$ are automata, they both only have one initial vertex. So let $V^i(\Gamma) = \{x\}$ and $V^i(\Lambda) = \{y\}$ and denote by $V' = (x, y)\cup \bar{V}(\Gamma)\times\bar{V}(\Lambda)$. Let $C$ be the constant from Theorem \ref{main_structure} such that 
\[
\abs{E(\Theta^s)}\leq C\cdot m^2\cdot n.
\]
Let $s_{\Gamma}:\mathbb{S}_{\Gamma}\to \Gamma$ and $s_{\Lambda}:\mathbb{S}_{\Lambda}\to \Lambda$ be the graphs from Theorem \ref{long_cycles_graph}. By Theorem \ref{long_cycles_algorithm} we may assume that we have precomputed these graphs. We may also assume that we have computed the subsets $V_{\mathbb{S}}$ and $E_{\mathbb{S}}$ from Proposition \ref{adjacency_partition} along with the partitions $V_{\mathbb{S}} = V_1 \sqcup ...\sqcup V_c$ and $E_{\mathbb{S}} = E_1\sqcup ...\sqcup E_c$.

Denote by $A\subset L(\Theta)$ the set of paths $g:I\to \Theta$ such that $t(g)\in V^f(\Theta^s)$ and by $B\subset L(\Theta)$ the set of paths $g:I\to \Theta$ such that $t(g)\in V^f(s(\mathbb{S}))$. Then, by definition of $\Theta^s$, we have $L(\Theta) = A\cup B$. Thus, determining if $L(\Theta)\neq \emptyset$ is equivalent to determining if $A\neq \emptyset$ or $B\neq\emptyset$.

The logical structure of the proof is as follows. Let $\Theta^{\circ}\subset \Gamma\times_{\Delta}\Lambda$ be the connected component containing the vertex $(x, y)$. We will first compute a graph $\Theta_k\subset \Gamma\times_{\Delta}\Lambda$ such that $\Theta_k$ contains the subgraph $\Theta^s\cap \Theta^{\circ}$ and satisfies $\abs{E(\Theta_k)}\in O(m^q\cdot n)$ for some $q\geq 0$. It will also be shown that $V(\Theta_k)\cap (V^f(\Gamma)\times V^f(\Lambda))\neq \emptyset$ if and only if $A\neq \emptyset$. Thus we will be able to decide if $A\neq \emptyset$ in $O(m^q\cdot n)$ time for some $q\geq 0$. Then, we will show that using this graph and Theorem \ref{unary_generalisation} we may also decide if $B\neq \emptyset$ in $O(m^q\cdot n)$ time for some $q\geq 0$.

Before showing how to compute $\Theta_k$ we will need a useful procedure. We will call the following algorithm the \emph{bounded forced path procedure}. Given as input a pair $(e, f)\in E(\Gamma)\times E(\Lambda)$ with $\gamma(e) = \lambda(f)$ and an integer $l$, do the following. Denote by $g:I\to \Gamma\times_{\Delta}\Lambda$ the path traversing only $(e, f)$. If $t(e, f)\in \bar{V}(\Gamma)\times\bar{V}(\Lambda)$ or $l = 1$, then we terminate and return the path $g$. If not, then $l\geq 2$ and $\odeg(e) = 1$ or $\odeg(f) =1$. Since $\gamma$ and $\lambda$ are forwards immersions, there is at most one element $(e', f') \in \fstar((t(e), t(f)), \Gamma\times_{\Delta}\Lambda)$. If no such element exists, then we also terminate and return the path $g$. If instead there is such a pair $(e', f')$, then let $g':I\to \Gamma\times_{\Delta}\Lambda$ be the output of the bounded forced path procedure with input $(e', f')$ and $l-1$. Return the path $g*g'$.

The following complexity bound follows from the fact that for all $u, v\in V(\Gamma)\times V(\Lambda)$, the set $\fstar((u, v), \Gamma\times_{\Delta}\Lambda)$ can be computed in $O(m)$ time.

\begin{lemma}
\label{forced_path_algorithm}
The bounded forced path procedure runs in $O(m\cdot l)$ time.

\end{lemma}

We are now ready to compute the graph $\Theta_k$ as promised. We will in fact do more than this, we will compute the following finite sequences:
\begin{enumerate}
\item graphs $\Theta_0, \Theta_1, ..., \Theta_k$ with $\Theta_i\subset \Gamma\times_{\Delta}\Lambda$ for all $i\geq 0$ 
\item sets of edges $\ndir_0, \ndir_1, ..., \ndir_k$, $\odir_0, \odir_1, ..., \odir_k$ with $$\ndir_i, \odir_i\subset \cup_{v\in V'}\fstar(v, \Gamma\times_{\Delta}\Lambda),$$ such that $\ndir_i\cap \odir_i = \emptyset$ and $\abs{\odir_i}<\abs{\odir_{i+1}}$ for all $i\geq 0$.
\end{enumerate}
By definition, we have $k\leq \abs{\bar{E}(\Gamma)}\cdot \abs{\bar{E}(\Lambda)}$ and so $k\in O(m^2)$.

We begin the computation by letting $\text{ndir}_0= \fstar((x, y), \Gamma\times_{\Delta}\Lambda)$, $\text{odir}_0=\emptyset$ and $\Theta_0 = \emptyset$. The set $\ndir_i$ will contain the new directions that we need to explore at step $i + 1$, meanwhile, the set $\odir_i$ will contain the old directions that we have already explored before step $i + 1$. So suppose now that we have computed $\Theta_{i-1}$, $\ndir_{i-1}$ and $\odir_{i-1}$. If $\ndir_{i-1} = \emptyset$, then we are done and $k = i-1$. If not, choose some element $(e, f)\in \ndir_{i-1}$. We have precomputed $G$, so we first check if 
\[
(e, f) \in d(E(G))\subset \cup_{v\in V'}\fstar(v, \Gamma\times_{\Delta}\Lambda). 
\]
We now have two cases to consider.

\textbf{Case 1:} $(e, f)\in E_{\mathbb{S}}$.

If $(e, f)\in E_{\mathbb{S}}$, then $(e, f)\notin E(\Theta^s)$ by definition. Now let $j$ be the integer such that $(e, f)\in E_j$. By Proposition \ref{adjacency_partition}, the vertices $V_j$ are precisely the vertices in $\bar{V}(\Gamma)\times \bar{V}(\Lambda)$ that can be reached by a path starting with $(e, f)$ in $s(\mathbb{S})$. So let:
\begin{align*}
\odir_i &= \odir_{i-1}\cup (e, f)\\
\ndir_i &= (\ndir_{i-1} - (e, f))\cup \bigcup_{v\in V_j}(\fstar(v, \Gamma\times_{\Delta}\Lambda) - (\fstar(v, \Gamma\times_{\Delta}\Lambda)\cap \odir_i))\\
\Theta_i &= \Theta_{i-1}.
\end{align*}

\textbf{Case 2:} $(e, f)\notin E_{\mathbb{S}}$.

Let $g$ be the output of the bounded forced path procedure with input the pair $(e, f)$ and the integer $C\cdot m^2\cdot n$. By Lemma \ref{forced_path_algorithm} this may be done in $O(m^3\cdot n)$ number of steps. Then we are left with two subcases.

\textbf{Subcase 2.1:} $t(g)\in \bar{V}(\Gamma)\times \bar{V}(\Lambda)$.

If this is the case, then we set:
\begin{align*}
\odir_i &= \odir_{i-1}\cup (e, f)\\
\ndir_i &= (\ndir_{i-1} - (e, f)) \cup (\fstar(t(g), \Gamma\times_{\Delta}\Lambda) - (\fstar(t(g), \Gamma\times_{\Delta}\Lambda)\cap \odir_i))\\
\Theta_i &= \Theta_{i-1}\cup g.
\end{align*}

\textbf{Subcase 2.2:} $t(g)\notin \bar{V}(\Gamma)\times\bar{V}(\Lambda)$.

Then let $j$ be the largest index such that $g(j)\in V^f(\Gamma)\times V^f(\Lambda)$. If no such integer exists, then set $j = 0$. We update our sets as follows:
\begin{align*}
\odir_i &= \odir_{i-1}\cup (e, f)\\
\ndir_i &= \ndir_{i-1} - (e, f)\\
\Theta_i &= \Theta_{i-1}\cup g[:j].
\end{align*}

We have the following properties of $\Theta_k$ and $\odir_k$ that follow from the construction:
\begin{enumerate}
\item\label{itm:theta_component} $\Theta^s\cap \Theta^{\circ}\subset \Theta_k$,
\item\label{itm:cycles} $E(s(\mathbb{S}))\cap E(\Theta_k) = \emptyset$,
\item\label{itm:accepting} for each $v\in \bar{V}(\Gamma)\times\bar{V}(\Lambda)$, there exists an accepting path $g:I\to \Theta$ traversing $v$ if and only if $v\in \bigcup_{(e, f)\in \odir_k} o(e, f)$.
\end{enumerate}

By \ref{itm:theta_component}, we have $(V^f(\Gamma)\times V^f(\Lambda))\cap V(\Theta_k)\neq \emptyset$ if and only if $A\neq\emptyset$. Hence, we have shown that may decide if $A\neq \emptyset$ in $O(m^q\cdot n)$ time for some $q\geq 0$. It remains to show that we can determine if $B\neq \emptyset$ in $O(m^q\cdot n)$ time for some $q\geq 0$.

Suppose that $B\neq \emptyset$. Then there is a pair of paths $g_1, g_2:I\to \Theta$ such that $o(g_1) = (x, y)$, $t(g_1) = o(g_2)\in \bar{V}(\Gamma)\times\bar{V}(\Lambda)$, $t(g_2)\in V^f(\Theta)$ and $g_2$ lifts to $s:\mathbb{S}\to \Theta$. By construction of \ref{itm:accepting}, it suffices to determine if there is a path $g:I\to \Theta$ such that $o(g)\in \bigcup_{(e, f)\in \odir_k}o(e, f)$, $t(g)\in V^f(\Theta)$ and $g$ lifts to $s:\mathbb{S}\to \Theta$. In order to determine if such a path exists, we essentially repeat the algorithm from Proposition \ref{long_graph_directions}. Let $(u, v)\in \bigcup_{(e, f)\in \odir_k}o(e, f)$ and $S_{\Gamma}\subset \mathbb{S}_{\Gamma}$, $S_{\Lambda}\subset \mathbb{S}_{\Lambda}$ a pair of connected components. Set
\begin{align*}
V^i(S_{\Gamma}) &= s_{\Gamma}^{-1}(u)\cap V(S_{\Gamma}),\\
V^i(S_{\Lambda}) &= s_{\Lambda}^{-1}(v)\cap V(S_{\Lambda}),\\
V^f(S_{\Gamma}) &= s_{\Gamma}^{-1}(V^f(\Gamma))\cap V(S_{\Gamma}),\\
V^f(S_{\Lambda}) &= s_{\Lambda}^{-1}(V^f(\Lambda))\cap V(S_{\Lambda}).
\end{align*}
Then $\core(S_{\Gamma}\times_{\Delta}S_{\Lambda})$ has an accepting path if and only if there is a path $g:I\to \Theta$ such that $o(g) = (u, v)$, $t(g)\in V^f(\Theta)$ and $g$ lifts to $s:\core(S_{\Gamma}\times_{\Delta}S_{\Lambda})\to \Theta$. Thus, by iterating over all such vertices and pairs of components, we may decide if $B\neq \emptyset$. Now, just as in the proof of Proposition \ref{long_graph_directions}, we may use Theorem \ref{unary_generalisation} to obtain the required bounds.
\end{proof}

\bibliographystyle{abbrv}
\bibliography{bibliography}

\end{document}